# Shape calculus for fitted and unfitted discretizations: domain transformations vs. boundary-face dilations

Martin Berggren

*Department of Computing Science, Umeå University, Sweden*

Dedicated to the memory of Roland Glowinski

May 25, 2023

**Abstract**

Shape calculus concerns the calculation of directional derivatives of some quantity of interest, typically expressed as an integral. This article introduces a type of shape calculus based on localized dilation of boundary faces through perturbations of a level-set function. The calculus is tailored for shape optimization problems where a partial differential equation is numerically solved using a fictitious-domain method. That is, the boundary of a domain is allowed to cut arbitrarily through a computational mesh, which is held fixed throughout the computations. Directional derivatives of a volume or surface integral using the new shape calculus yield purely boundary-supported expressions, and the involved integrands are only required to be element-wise smooth. However, due to this low regularity, only one-sided differentiability can be guaranteed in general. The dilation concept introduced here differs from the standard approach to shape calculus, which is based on domain transformations. The use of domain transformations is closely linked the the use of traditional body-fitted discretization approaches, where the computational mesh is deformed to conform to the changing domain shape. The directional derivatives coming out of a shape calculus using deforming meshes under domain transformations are different then the ones from the boundary-dilation approach using fixed meshes; the former are not purely boundary supported but contain information also from the interior.

## 1 Introduction

Gradient-based shape optimization is a powerful tool for engineering design, as it enables algorithmic exploration of larger design spaces than is possible by manual parameter studies. Gradient-based techniques are particularly effective for cases when there is only a small number of quantities of interest but a very large number of design variables, since by employing the adjoint-variable technique, the calculation time will essentially be independent of the number of design variables for any fixed discretization level.

When considering gradient-based shape optimization for a particular problem, a central part is the *shape calculus* (or the *shape sensitivity analysis*), the calculation of directional derivatives of the quantities of interest — objective functions and constraints — with respect to the design variables. These directional derivatives are then used to define the gradients needed by the optimization algorithm. Conceptually, shape calculus relies on a calculus of variation with respect to boundary modifications, and methods on how to carry out such calculations are very well established by now; standard references are the extensive monographs by Delfour & Zolésio [9] and Sokolowski & Zolésio [20].

As his graduate student in the early 1990s, Roland Glowinski taught me how to carry out sensitivity analysis for various control problem associated with partial differential equations. He pointed out the need for *accuracy* in such calculations. In particular, he strongly promoted the idea of carrying out the calculations from scratch in the *fully discrete case*, so that the final expression will be the exact directional derivative, up to round-off, of the discrete objective function; the argument is that otherwise it will be impossible to guarantee convergence of the iterations in optimization algorithms of descent type. This approach is known as the "discrete approach," or "discretize-then-differentiate." Glowinski & He [12] provide an example of non-convergence when such a scheme is not followed.

My experience in carrying out shape optimization, in particular for wave-propagation problems [2,22] (e.g.), supports the benefits of the fully discrete shape calculus approach to achieve a robust performance of the optimization algorithm. However, when a mesh deformation algorithm is used to modify the location of interior mesh nodes with respect to a reshaped boundary, the dependency on this algorithm is a complication with the fully discrete approach. The exact directional derivative will depend on the "mesh sensitivities," that is, on the mapping between the displacements of the boundary and internal nodes as well as on the Jacobian transpose of this mapping [3, § 3.2]. This complication introduces at the very least a software dependency between tools that may be more convenient to keep separated.

A simpler approach from an implementation point of view is to use directional-derivative expressions involving only boundary quantities, coming from an analysis of the problem before discretization. This approach is known as the "continuous approach" or "differentiate-then-discretize" and ignores the dependency on the mesh deformation algorithm. In spite of my own preferences, I have to admit to have accepted this approach in a full 3D case to simplify the implementation [18]; we obtained good results in this case, but I believe that a fully discrete approach is still to be preferred, if possible.

Backing away from these well-established concepts in order to take a fresh look at the subject; what is the basic reason for the occurrence of these varying directional-derivative expressions? I believe it goes all the way back to the decision *how to treat geometry as a variable* with respect to which differentiation is carried out. The standard, universally used approach, covered for instance in the above mentioned monographs [9, 20], is build on *domain transformations*, analogous to how the mechanics of continuous media is handled. It is so universally accepted that it is easy to forget that the use of transformations is a choice! Using domain transformations, an integral over a modified geometry can be transformed back to the unmodified geometry, which means that the parameters of the transformation will appear in the integrand. Then, ordinary differential calculus can be applied to the transformed integrand.

The classical approach to shape calculus using domain transformation is particularly well suited to traditional numerical methods using body-fitted meshes, since domain transformations are also useful to modify the interior mesh in order to preserve mesh quality. However, mesh deformation algorithms are an annoyance in themselves. The algorithms can have robustness issues, particularly in three space dimensions, and large changes in the domain shape will generally require remeshing, adding an extra layer of complexity. For shape optimization, it is therefore of interest to consider *fictitious-domain methods*, in which the domain $\Omega$ of interest is embedded in a fixed hold-all domain $D$. No mesh deformations or remeshings are then needed; instead of fitting the mesh to the domain, the fictitious-domain methods instead restrict the computations to $\Omega$ and assign boundary conditions on $\partial\Omega$, which will be an interior boundary in $D$. Roland Glowinski was a noted contributor to the development of fictitious-domain methods [13, 14] (e.g.), and there have been a recent reemergent interest in a class of fictitious-domain methods known as CutFEM or XFEM [7]. We have demonstrated the effectiveness of the CutFEM for acoustic shape optimization [4, 5], and we also discovered some striking differences with respect to the use of mesh-deformation techniques, to be expanded on in this article. A first observation is that the domain-transformation approach to the treatment of geometry changes is *not* natural for fictitious-domain methods, since domain changes are then effectuated just by moving the boundary over a fixed background mesh. It is possible to artificially introduce a domain transformation, which is what we did [4, 5], but since the mesh is fixed, the transformation should be localized just around the mesh elements crossing the boundary in order to reflect what really happens in the computations.

The purpose of the current contribution is to introduce a shape calculus approach that is conceptually different from the domain-transformation idea, namely an approach built on *localized dilations of boundary faces*. The idea, introduced in § 5.2, constitutes a generalization of what is here called *Delfour dilations*. Delfour [10] introduced the concept of domain perturbations generated by dilations of lower-dimensional objects into Euclidean space. For instance, the dilation of a point can be used to generate a domain perforated by a ball, and the dilation of a curve to generate a domain perforated by a curved rod. Of interest here is dilations of the boundary of a bounded domain. As demonstrated by Delfour, this dilation can be used to define a shape derivative with respect to a uniform extension or contraction of the domain. Here we generalize this approach and apply it in the discretized case to localized dilations effectuated by a perturbation of a level-set function. The use of a level-set function to specify the geometry is a standard and attractive approach in shape optimization, which is why the current contribution could serve as a tool box for shape calculus in the context of fictitious-domain methods using level-set geometry descriptions.

In a recent article, Gangl and Gfrerer [11] also extend the Delfour dilations using, as here, perturbations of a level-set function in the discretized case. Their aim is to unify the topological and shape derivative into an object they call the *topological–shape derivative*. Gangl and Gfrerer's contribution is complementary to what is presented here. A main focus here is on the general directional-derivative expressions in Theorems 6.7 and 6.11, valid in both 2 and 3 dimensions, with minimal regularity requirements on the integrands. In contrast, Gangl and Gfrerer's main focus is on the interesting concept of defining a combined topological–shape derivative, which they apply to a 2-dimensional tracking-type model problem in order to work out detailed derivative expressions on the matrix–vector level.

To highlight the contrast of the approach introduced here to the traditional approach, the standard shape calculus using domain transformations is reviewed in § 4. The analysis is applied to the simple model shape optimization problem introduced in § 2. In particular, we highlight the different directional derivative expressions that are obtained for the problem before and after discretization when deforming meshes are used. In contrast, § 5.4 shows that the fully discrete shape derivative is much simpler when a fictitious-domain method with a fixed mesh is used, and the final expression agrees more closely with the simple boundary expression for the directional derivative before discretization.

This article constitutes a kind of continuation of contribution [3], published in the proceedings of the conference arranged on the occasion of Roland Glowinski's 70th birthday. The previous article presented a unified shape calculus for the cases before and after discretization when domain transformations are employed. It showed that the exact directional derivative in the discrete case when mesh deformations are used cannot avoid the changes in mesh and solutions taking place inside the domain. It also highlighted which terms that will be ignored when using the simpler boundary integral expression coming from the "differentiate-than-discretize" approach. In contrast, the present contribution shows that the directional derivatives when using fictitious-domain methods and fixed meshes can be expressed as boundary expressions, a result that even more increases the attractiveness of such methods for shape optimization.

The proofs of the shape calculus formulas for surface dilations, although conceptually not difficult, are unfortunately quite long and somewhat tedious. To ease the exposition for readers more interested in the big picture and the resulting formulas, these proofs are collected in a separate § 6.

## 2  A model shape optimization problem

As a model problem to illustrate the differences between the shape calculus for various cases, we consider *compliance minimization* for the following Poisson problem,

$$
\begin{aligned}
-\Delta u &= r && \text{in } \Omega, \\
u &= 0 && \text{on } \Gamma_{\text{D}}, \\
\alpha u + \frac{\partial u}{\partial n} &= 0 && \text{on } \Gamma_{\text{R}},
\end{aligned}
\tag{2.1}
$$

where $\Omega$ is an open, bounded, and connected domain in $\mathbb{R}^d$ ($d = 2$ or 3), $r \in H^1(\mathbb{R}^d)$ is a given function, $\alpha \geq 0$ a constant, and the domain boundary $\partial\Omega$ comprises the nonoverlapping parts $\Gamma_{\text{D}}$ and $\Gamma_{\text{R}}$, both of positive $d - 1$-dimensional measure.

The objective function is

$$
J(\Omega) = \int_\Omega r u \, dV,
\tag{2.2}
$$

and the optimization problem is to find the shape of the domain $\Omega$ that minimizes $J$, typically subject to constraints such as a bound on the volume of $\Omega$. Well-posedness of this problem requires some regularization scheme such as a suitable restriction on the class of admissible domains. However, here our interest will only be the shape calculus with respect to changes in $\Omega$ effectuated either by domain transformations (§ 4) or by dilations (§ 5).

*Remark* 2.1. The name "compliance minimization" is borrowed from structural mechanics. In the corresponding problem for linear elasticity, the objective function corresponds to the *mechanical compliance* of the structure, or, equivalently, to the *work* carried out on the structure by the external forces $r$.

## 3  Shape functions and shape calculus using domain paths

We define the *hold-all* $D \subset \mathbb{R}^d$ as the point set within which all feasible domains will be contained. For the purpose of this article, it will be convenient to constrain $D$ to be nonempty, open, bounded, and simply connected. A *shape function* is a function $J: \mathcal{U} \to \mathbb{R}$, defined on a family $\mathcal{U}$ of admissible subsets compactly embedded in $D$. Typical example shape functions are the domain and boundary integrals

$$J_1(\Omega) = \int_\Omega f \, dV, \qquad J_2(\Omega) = \int_{\partial\Omega} f \, dS, \tag{3.1}$$

where $f \in L^1(D)$ and $f \in W^{1,1}(D)$ for $J_1$ and $J_2$, respectively; further requirements are specified below when needed. To carry out shape calculus on such functions, we will consider *domain paths* $t \to \Omega_t$ in $\mathcal{U}$ and, in particular, limits of the type

$$\lim_{t \to 0^+} \frac{J_k(\Omega_t) - J_k(\Omega)}{t}. \tag{3.2}$$

Depending on how these paths are constructed, the aim of the shape calculus is to identify such limits as directional derivatives in order to supply gradient information to the optimization algorithm.

As reviewed in § 4, the classical way of constructing such paths relies on domain transformations. An alternative, based on boundary dilations is presented in § 5. In particular, § 5.2 introduces a dilation method suitable for fictitious-domain methods.

## 4 Domain paths using transformations

Here, in order to construct domain paths, we fix a given reference domain $\Omega \subset\subset D$ with a Lipschitz boundary and generate the path through a family of smooth homeomorphisms $T_t: D \to \mathbb{R}^d$, parametrized by $t \geq 0$. In this way, shape calculus can be converted to ordinary differential calculus by a change of variables based on $T_t$.

Specifically, we consider objective functions (3.1) defined on the mapped domain $\Omega_t = T_t(\Omega)$ and with integrands given by a family of functions $t \mapsto f_t$, where each $f_t \in L^1(D)$ for $J_1$ and each $f_t \in W^{1,1}(D)$ for $J_2$. A change of variables yields that

$$J_1(\Omega_t) = \int_{\Omega_t} f_t \, dV = \int_\Omega f_t \circ T_t |\det DT_t| \, dV, \tag{4.1a}$$

$$J_2(\Omega_t) = \int_{\partial\Omega_t} f_t \, dS = \int_{\partial\Omega} f_t \circ T_t |(DT_t)^{-T} \boldsymbol{n}| |\det DT_t| \, dS, \tag{4.1b}$$

where $DT_t$ is the Jacobian matrix of the transformation and $\boldsymbol{n}$ the normal field on the boundary. The transformation now appears explicitly in the integrands, which means that the integrals can be differentiated with respect to $t$ using ordinary calculus.

Now consider the special case of a transformation such that $T_t(\Omega) = \Omega$ for each $t \geq 0$; note that such a transformation is not necessarily just the identity! When the integrand function is fixed, that is, there is an $f$ such that $f_t = f$ for each $t$, the shape functions will be unchanged under such a transformation. That is, a transformation that just involves a rearrangement of the points inside the domain will not change the value of such shape functions; in this case, difference quotient $(J(\Omega_t) - J(\Omega))/t$ will therefore vanish. Indeed, this property motivates the terminology *shape* functions for integrals (3.1); the only transformations that may yield nontrivial derivatives are those that modify $\partial\Omega$. This transformation property is also reflected in the so-called *structure theorem* [9, Thm. 3.6, Ch. 9], which says that the *gradients* of the shape functions (3.1) can be identified as distributions on $\mathbb{R}^d$ with their support contained in $\partial\Omega$.

Now, how can the transformations $T_t$ be defined? One option is employ an explicit parametrization. For instance, in the case $D$ is a cube ($d = 3$) or a square ($d = 2$), transformations of a reference domain $\Omega \subset D$ can be effectuated through *free-form deformations*, a technique originating in computational geometry [16, 19]. Here a grid of control points are distributed throughout $D$, and the parameters of the transformation will be displacements in each coordinate direction of these control points. (The directional derivatives will be computed with respect to these displacements.) To utilize formulas (4.1), the point displacements need to be extended to a transformation on $D$, for instance using Bernstein polynomials or non-uniform rational B-splines (NURBS), depending on whether the control points are uniformly distributed in $D$ or not.

A more common and more flexible approach is to define each transformation using a so-called *velocity field* $\boldsymbol{V}: D \to \mathbb{R}^d$. The simplest use of the velocity field is as a *perturbation of identity*; that is, for each

$x \in \Omega$,
$$T_t(x) = x + tV(x). \tag{4.2}$$

Another option is to use the so-called *velocity method*, in which the evaluation of the transformation at each $x \in \Omega$ is the flow of the vector field initiated at $x$. That is, $T_t(x) = X_x(t)$, where $t \mapsto X_x(t)$ satisfies

$$\begin{aligned} \frac{dX_x}{dt} &= V \circ X_x \qquad t > 0, \\ X_x(0) &= x. \end{aligned} \tag{4.3}$$

(It is straightforward to generalize equation (4.3) to use a nonautonomous vector field $(x, t) \mapsto V(x, t)$, a generality not needed here.) The transformations defined through perturbations of identity and the velocity method agree to first order at $t = 0$ and will thus provide the same first derivatives at $t = 0$ of shape functions (4.1).

Both integrands in the transformed objective functions (4.1) involve the compound function $f_t \circ T_t$ as well as a geometric term, $|\det DT_t|$ and $|(DT_t)^{-T} n||\det DT_t|$, respectively. A differentiation of the compound function with respect to $t$ yields the *material derivative*

$$\dot{f}(x) \stackrel{\text{def}}{=} \lim_{t \to 0^+} \frac{(f_t \circ T_t)(x) - f(x)}{t} = \frac{d^+}{dt}\bigg|_{t=0} (f_t \circ T_t)(x). \tag{4.4}$$

The differentiated geometric terms can, after a few steps of algebra, be shown to satisfy

$$\begin{aligned} \frac{d^+}{dt}\bigg|_{t=0} \det DT_t &= \nabla \cdot V, \\ \frac{d^+}{dt}\bigg|_{t=0} \det DT_t |(DT_t)^{-T} n| &= \nabla \cdot V - \frac{\partial}{\partial n}(n \cdot V) \stackrel{\text{def}}{=} \nabla_{\mathrm{T}} \cdot V; \end{aligned} \tag{4.5}$$

the last expression is the *tangential divergence*[1] of the velocity field along $\partial\Omega$. Using formulas (4.4) and (4.5) together with the product rule, a differentiation under the integral signs of objective functions (4.1) yields the directional derivatives

$$dJ_1(\Omega; V) = \int_\Omega \left(\dot{f} + f \nabla \cdot V\right) dV, \tag{4.6a}$$

$$dJ_2(\Omega; V) = \int_{\partial\Omega} \left(\dot{f} + f \nabla_{\mathrm{T}} \cdot V\right) dS, \tag{4.6b}$$

using the notation $f = f_t|_{t=0}$, a convention that will be followed henceforth.

As we have seen, formulas (4.6) are rather straightforward consequences of defining domain paths using transformation (4.2) and hold as long as all constituents are well-defined and integrable. Formula (4.6a) holds for $V \in W^{1,\infty}(D)^d$ and when the compound function $t \mapsto f_t \circ T_t$ is continuous with values in $L^1(\Omega)$ and differentiable with values in $L^1(\Omega)$, and formula (4.6b) for $V \in W^{2,\infty}(D)^d$ and when $t \mapsto f_t \circ T_t|_{\partial\Omega}$ is continuous with values in $L^1(\partial\Omega)$ and differentiable with values in $L^1(\partial\Omega)$

Under different and stronger conditions, it is possible to devise a different set of formulas for the directional derivatives. These formulas, proved by Delfour & Zolésio [9, Ch. 9, §4], contain the *shape derivative* $f'$, that is, the derivative of the map $t \mapsto f_t$ evaluated at $t = 0$, and the local summed curvature $\kappa = \nabla_{\mathrm{T}} \cdot n$ of $\partial\Omega$,

$$dJ_1(\Omega; V) = \int_\Omega f' \, dV + \int_{\partial\Omega} f V \cdot n \, dS, \tag{4.7a}$$

$$dJ_2(\Omega; V) = \int_{\partial\Omega} f' \, dS + \int_{\partial\Omega} \left(\frac{\partial f}{\partial n} + \kappa f\right) V \cdot n \, dS. \tag{4.7b}$$

The proofs by Delfour & Zolésio require $V \in C^1(D)$ and the following additional assumptions.

---

[1] The definition can also be written $\nabla_{\mathrm{T}} \cdot V = \nabla \cdot V - n \cdot (\nabla V) n = \nabla \cdot V - n \cdot ((n \cdot \nabla)V)$, where $\nabla V$ is the Jacobian of $V$.

- For formula (4.7a): the function $t \mapsto f_t$ is continuous with values in $W^{1,1}(D)$ and differentiable with values in $L^1(D)$.

- For formula (4.7b): $\Omega$ is of class $C^2$ and the function $t \mapsto f_t$ is continuously differentiable with values in $H^2(D)$.

Both sets of formulas, the "weak" directional derivatives (4.6) as well as the "strong" (4.7), are useful for shape optimization problems. For the problem introduced in § 2, the analysis in § 4.1 below will assume sufficient regularity for expressions (4.7) tho hold. However, we will see that when applying a standard finite-element discretization approach using deforming meshes, the regularity assumptions fail to hold. As will be discussed in § 4.2, formulas (4.6) are the appropriate choice in the discrete case when using deforming meshes. Note that these considerations are contingent on the use of domain transformations. Indeed, we will see that when creating domain paths using dilations, the role of these formulas will in a sense be reversed, so that a version of formulas (4.7) will then be the appropriate choice in the discrete case!

### 4.1 Shape calculus for the model problem before discretization

We consider the problem of § 2 and a domain path $t \mapsto \Omega_t$ generated by domain transformation (4.2) associated with a smooth velocity field $V$ vanishing on $\Gamma_\mathrm{D}$. The standard variational formulation of the Poisson problem (2.1) on the deformed domain $\Omega_t$ is

$$u(t) \in W(t) \text{ such that}$$
$$\int_{\Omega_t} \nabla v \cdot \nabla u(t)\, dV + \alpha \int_{\Gamma_{\mathrm{R},t}} v u(t)\, dS = \int_{\Omega_t} rv\, dV \qquad \forall v \in W(t), \tag{4.8}$$

where $W(t)$ is the subset of functions in $H^1(\Omega_t)$ with vanishing trace on $\Gamma_\mathrm{D}$. Note that $\Gamma_\mathrm{D}$ is fixed under this domain transformation, since we have chosen $V$ to vanish on $\Gamma_\mathrm{D}$, whereas $\Gamma_\mathrm{R}$ may be displaced by the perturbation. Objective function (2.2) evaluated on the deformed domain will be

$$J(\Omega_t) = \int_{\Omega_t} ru(t)\, dV. \tag{4.9}$$

The conventional shape calculus approach for this kind of problems is to view objective function (4.9) as a composite mapping involving the intermediate *state variable* $u(t)$, which will be differentiated with respect to $t$. One way to carry out these calculations is to apply formulas (4.6), which will involve the material derivative $\dot u$. The existence of $\dot u$ can be rigorously established, as shown by Sokolowski & Zolézio [20, § 2.28–2.29] for similar problems. Here, however, we will demonstrate the use of formulas (4.7). Since these involve boundary traces of derivatives, we need to assume enough smoothness of $\Omega$ for full elliptic regularity of problem (4.8) to hold; that is, we assume that $u(t) \in H^2(\Omega_t) \cap W(t)$. Note that this assumption may not be completely justified, since the presence of the mixed boundary conditions in problem (2.1) can introduce a local loss of regularity at the interface between $\Gamma_\mathrm{D}$ and $\Gamma_\mathrm{R}$. Moreover, our calculations will involve the shape derivative $u' = \lim_{t \to 0+} \bigl(u(t) - u(0)\bigr)/t$, whose existence, properties, and the regularity requirements necessary for a rigorous analysis are nontrivial and delicate. Thus, our calculation will be merely formal.

These complications, together with the fact that the final directional-derivative expression will not depend on a differentiated state variable, has lead several authors to search for ways to rigorously establish a shape calculus without requiring a differentiation of the state. Using a Lagrangian formalism, Delfour & Zolésio [9, Ch. 10] discuss ways of carrying out such a task. A more general framework for shape calculus that circumvents the introduction of $\dot u$ or $u'$ has been introduced by Ito, Kunish, and Peichl [15]. Based on the introduction of a so-called averaged adjoint equation, Sturm [21] has devised another state-differentiation-free approach with relaxed continuity assumptions. Finally, in a recent contribution, Laurain et al. [17] discuss, compare, and systematize several Lagrangian-based methods bypassing the need for state differentiation.

Since the objective here, however, is just to illustrate different shape calculus approaches and compare the final expressions, we will follow a quick, formal approach and simply assume that $u'$ exists and is regular enough for all required calculations to be valid. In particular, in this section we assume that $u' \in W$.

Under this assumption, the directional derivative of objective function (4.9) at $t = 0$ can be computed by first observing that the solution to problem (4.8) satisfies

$$\int_{\Omega_t} |\nabla u(t)|^2 \, dV + \alpha \int_{\Gamma_{R,t}} u(t)^2 \, dS = \int_{\Omega_t} ru(t) \, dV, \tag{4.10}$$

which means that objective function (4.9) at $t = 0$ can be written

$$J(\Omega) = \int_\Omega ru \, dV = 2\int_\Omega ru \, dV - \int_\Omega |\nabla u|^2 \, dV - \alpha \int_{\Gamma_R} u^2 \, dS. \tag{4.11}$$

Formulas (4.7) and the product rule then yield that

$$\begin{aligned}
dJ(\Omega; V) &= 2\int_\Omega ru' \, dV - 2\int_\Omega \nabla u' \cdot \nabla u \, dV + \int_{\partial\Omega} \boldsymbol{n} \cdot \boldsymbol{V} \left(2ru - |\nabla u|^2\right) dS \\
&\quad - 2\alpha \int_{\Gamma_R} u'u \, dS - \alpha \int_{\Gamma_R} \boldsymbol{n} \cdot \boldsymbol{V} \left(\frac{\partial}{\partial n} u^2 + \kappa u^2\right) dS \\
&= \int_{\Gamma_R} \boldsymbol{n} \cdot \boldsymbol{V} \left(2ru - |\nabla u|^2 - \alpha \left(\frac{\partial}{\partial n} u^2 + \kappa u^2\right)\right) dS,
\end{aligned} \tag{4.12}$$

where the vanishing of $V$ on $\Gamma_D$ and variational form (4.8) with $v = u'$ has been used in the second equality. Note that we have also used that $r' = 0$ since $r \in H^1(D)$ does not depend on $t$.

Directional derivative (4.12) vanishes whenever $\boldsymbol{n} \cdot \boldsymbol{V} = 0$, which shows that objective function (4.9) in this case indeed is a shape function. However, as we will see, the situation changes after finite-element discretization.

*Remark* 4.1. Note that the directional derivative (4.12) only depends on the solution $u$ and the velocity field $V$. This is a rather special case for problems having certain symmetries. Other choices of objective functions will necessitate the introduction of an *adjoint* state as well. In the current case, the state equation does double duty also as adjoint equation.

### 4.2 Discrete shape calculus with body-fitted meshes and mesh deformations

Standard finite-element methods use meshes that are fitted to the domain, at least approximately. Typically, the mesh boundary points are located exactly on the domain boundary, whereas the union of the mesh boundary surfaces will in general only approximate the domain boundary if it is curved. We introduce such a nondegenerate triangulation $\mathscr{T}_h$ using simplicial elements. The union of the elements in $\mathscr{T}_h$ generates a approximating computational domain $\Omega_h$ containing boundaries $\Gamma_{D,h}$ and $\Gamma_{R,h}$. Parameter $h > 0$ specifies the largest diameter of any element in $\mathscr{T}_h$. Moreover, by $\mathscr{S}_h$ and $\mathscr{E}_h$ we denote the set of faces and subfaces of the elements in $\mathscr{T}_h$. For simplicity of exposition, in this section we will specialize to 3 space dimensions, $d = 3$. Thus, the elements, faces, and subfaces are here tetrahedrons, triangles, and line segments. The changes needed to interpret the formulas below to $d = 2$ are rather straightforward. Based on the triangulation, we introduce finite-element function spaces $W_{h,p}$ of continuous functions that are polynomials of at most degree $p \geq 1$ on each element $K \in \mathscr{T}_h$ and that vanish on $\Gamma_{D,h}$.

The use of a domain-fitted mesh means that it has to be regenerated or adjusted when the shape of the domain changes. Gradient-based optimization algorithms assume smoothness in the evolution of the objective function under changes in the decision variables, which is why a complete regeneration of the mesh is usually performed only when necessary, as it will introduce a "meshing noise" in the evolution of the objective function. However, the domain transformation idea, which above was introduced just as a tool to generate domain paths for the definition of shape derivatives, fits very well also with the need to modify the mesh under shape changes. We thus apply transformation (4.2) with $V = V_h$, where $V_h \in W_{h,1}^d$ is the space of continuous, piecewise-linear vector-valued functions on the triangulated domain. This choice of space for $V_h$ assures that each face in $\mathscr{S}_h$ will remain planar under the transformation and, as long as $t$ is small enough, that the transformed mesh points still generate a valid mesh. Note that, similarly as in the case before discretization, $V_h$ vanishes on $\Gamma_{D,h}$.

The transformation deforms the mesh to generate a domain $\Omega_{h,t}$, and we define on the deformed mesh the deformed finite element space $W_{h,p}(t)$ and arrive at the following finite-element approximation

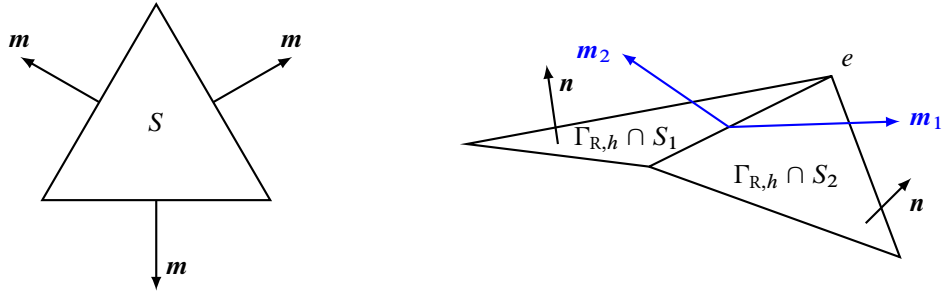

FIGURE 1. Left: the co-normal $\boldsymbol{m}$ is defined on each edge of a mesh face $S \in \mathscr{S}_h$ and is orthogonal to the face normal and the edge. Right: the co-normal jump over an edge $e \in \mathscr{E}_h$ interfacing two neighboring faces $S_1$, $S_2$ on a piece of $\Gamma_{R,h}$ is $[\![\boldsymbol{m}]\!] = \boldsymbol{m}_1 + \boldsymbol{m}_2$.

of variational problem (4.8):

$$u_h(t) \in W_{h,p}(t) \text{ such that}$$
$$\int_{\Omega_{h,t}} \nabla v_h \cdot \nabla u_h(t)\, dV + \alpha \int_{\Gamma_{R,h,t}} v_h u_h(t)\, dS = \int_{\Omega_{h,t}} r v_h\, dV \qquad \forall v_h \in W_{h,p}(t). \tag{4.13}$$

The discrete objective function is, consequently,

$$J_h(\Omega_{h,t}) = \int_{\Omega_{h,t}} r u_h(t)\, dV. \tag{4.14}$$

Now assume that we choose a $V_h \neq \mathbf{0}$ that vanishes on $\partial\Omega_h$. Corresponding transformation keeps the boundary nodes fixed but offsets the locations of internal nodes of the triangulation. Such a transformation will change the finite element solution $u_h$ (even though its order of approximation with respect to the exact solution $u$ is unchanged). Thus, $J_h$ can change, albeit slightly, with transformations that vanish on $\partial\Omega_h$, and it follows that $J_h$, as opposed to the $J$ in expression (4.11), *is not a strict shape function*. Indeed, as shown below in § 4.2.1, the directional derivative of objective function (4.14) at $t=0$ can be written in the form

$$\begin{aligned}dJ_h(\Omega_h; V_h) =& \\= & \int_{\Gamma_{R,h}} \boldsymbol{n} \cdot V_h \Big( 2ru_h - |\nabla u_h|^2 - \alpha\Big(\frac{\partial}{\partial n} u_h^2 + \kappa_h u_h^2\Big) \Big)\, dS - \sum_{e \in \mathscr{E}_h} \alpha \int_{e \cap \Gamma_{R,h}} [\![\boldsymbol{m}]\!] \cdot V_h u_h^2\, dl \\ & - \sum_{S \in \mathscr{S}_h^0} \int_S V_h \cdot [\![\boldsymbol{n}|\nabla u_h|^2]\!]\, dS \\ & + 2 \sum_{K \in \mathscr{T}_h} \int_K \Big( \nabla u_h \cdot \nabla(V_h \cdot \nabla) u_h - r V_h \cdot \nabla u_h \Big)\, dV + 2\alpha \int_{\Gamma_{R,h}} u_h (V_h \cdot \nabla) u_h\, dS.\end{aligned} \tag{4.15}$$

Before deciphering the notation, we note that the first integral in expression (4.15) looks like corresponding formula (4.12) before discretization. The other terms are due to the discretization. Note, in particular, that integrals 3 and 4 depend on $V_h$ and $u_h$ in the interior, reflecting the fact that $J_h$ is not a pure shape function.

The function $\kappa_h$ in the first integral of expression (4.15) is defined such that, for each $S \in \mathscr{S}_h$, $\kappa_h|_S = \nabla_T \cdot \boldsymbol{n}|_S$ is the summed curvature of $S$. Since we here assume planar mesh faces, $\kappa_h$ vanishes almost everywhere on $\Gamma_{R,h}$, so associated term in expression (4.15) could be removed. The term is kept in the expression in order to highlight the parallel to the case before discretization. Moreover, as we will see below, $\kappa_h$ will be nonzero if elements with curved faces at $\Gamma_{R,h}$ are used.

Discrete curvature contributions will come from the second integral in expression (4.15). This term constitutes a sum over the mesh edges on $\Gamma_{R,h}$ and contains the co-normal field $\boldsymbol{m}$, defined for each $S \in \mathscr{S}_h$ on the boundary $\partial S$ as indicated to the left in figure 1. As illustrated to the right in figure 1, on an edge $e$ interfacing two faces $S_1$, $S_2$ on $\Gamma_{R,h}$, we define the jump term

$$[\![\boldsymbol{m}]\!] = \boldsymbol{m}_1 + \boldsymbol{m}_2. \tag{4.16}$$

The integral of this jump term yields a discrete curvature measure; note that it vanishes if $S_1$ and $S_2$ are coplanar.

The third integral in expression (4.15) is taken over over all mesh faces $\mathscr{S}_h^0 \subset \mathscr{S}_h$ interior to the domain. Each such $S \in \mathscr{S}_h$ is a common face of two elements $K_1, K_2 \in \mathscr{T}_h$, and $[\![\boldsymbol{n}|\nabla u_h|^2]\!]$ denotes the jump of $|\nabla u_h|^2$ over the surface, as defined below in expression (4.25).

Finally, the last two integrals in expression (4.15) contain element-wise evaluations of the residual of the state equation (4.13) at $t = 0$, evaluations in which $V_h \cdot \nabla u_h \notin W_{h,p}$ replaces the test function. Assuming full elliptic regularity of weak solutions to the Poisson problem (2.1) and the convergence properties of finite-element solutions for properly chosen meshes, we expect integrals 3–5 of expression (4.15) to vanish as $h \to 0$.

### 4.2.1 How to obtain directional derivative (4.15)

To show expression (4.15), we use the same approach as in a previous contribution [3], based on a systematic use of the concept of material derivative, introduced already in definition (4.4). The chain rule yields the following relation of the material derivative and the shape derivative $f'$ of expressions (4.7),

$$f' = \dot{f} - \boldsymbol{V} \cdot \nabla f. \tag{4.17}$$

(The shape and material derivatives correspond to the *partial* and *total* derivatives in the mechanics of continuum media.)

We summarize here the rules for the material derivative that will be used to establish formula (4.15), and we refer to previous publication [3] for details. The restriction on $\Omega_t$ of a fixed $r \in H^1(D)$, such as the one in equation (4.8), has a vanishing shape derivative, $r' = 0$ (since it is fixed and not dependent on $t$). Thus, by relation (4.17),

$$\dot{r} = \boldsymbol{V} \cdot \nabla r. \tag{4.18}$$

The product rule holds for the material derivative, so for functions $g_1, g_2 : \Omega \to \mathbb{R}$,

$$(g_1 g_2)^{\cdot} = \dot{g}_1 g_2 + g_1 \dot{g}_2, \tag{4.19}$$

provided that the derivatives exist. We will need a product rule for $|\nabla g|^2$. The shape derivative, being a partial derivative, commutes with the gradient operator, but the material derivatives does not. Instead the rule

$$\left(|\nabla g|^2\right)^{\cdot} = 2\nabla \dot{g} \cdot \nabla g - 2\nabla g \cdot (\nabla \boldsymbol{V}) \nabla g \tag{4.20}$$

holds for $g \in H^1(\Omega)$ such that $\dot{g} \in H^1(\Omega)$ [3, eq. (30)]. (In Cartesian components, $\nabla \boldsymbol{V}$ is the Jacobian matrix of $\boldsymbol{V}$.)

Since the mesh is deformed, the support of each finite-element basis function $w$ for the functions in $W_{h,p}(t)$ will be transformed with the velocity field $\boldsymbol{V}_h$, which implies that $\dot{w} = 0$. Consequently, the material derivative of any finite-element function conforms to the finite-element space; that is, if $v_h \in W_{h,p}$ then $\dot{v}_h \in W_{h,p}$ [3, § 4, example 3]. However, we see from expression (4.17) that the shape derivative of a finite-element function does *not* conform. Indeed, $v_h' \notin W_{h,p}$ in general, since $\boldsymbol{V}_h \cdot \nabla v_h$ will generally not be continuous. This property is a source to the differences with respect to the case in § 4.1; recall that the property $u' \in W$ was used in expression (4.12) to arrive to the final expression!

We now possess the tools to establish formula (4.15). Proceeding similarly as in the case before discretization, we start by noticing that the solution to equation (4.13) satisfies

$$\int_{\Omega_{h,t}} |\nabla u_h(t)|^2 \, dV + \alpha \int_{\Gamma_{\mathrm{R},h,t}} u_h(t)^2 \, dS = \int_{\Omega_{h,t}} r u_h(t) \, dV, \tag{4.21}$$

which means that objective function (4.14) again can be written in the form

$$J_h(\Omega_h) = 2 \int_{\Omega_h} r u_h \, dV - \int_{\Omega_h} |\nabla u_h|^2 \, dV - \alpha \int_{\Gamma_{\mathrm{R},h}} u_h^2 \, dS. \tag{4.22}$$

Assuming that $u_h$ is material differentiable and recalling that the material derivative of the finite-element functions stays in the space, we utilize formulas (4.6) on objective function (4.22), which yields

$$\begin{aligned} dJ_h(\Omega_h; V_h) &= 2\int_{\Omega_h} (\dot{r} u_h + r\dot{u}_h + r u_h \nabla \cdot V_h) \, dV - \int_{\Omega_h} \left( (|\nabla u_h|^2)^{\cdot} + |\nabla u_h|^2 \nabla \cdot V_h \right) dV \\ &\quad - \alpha \int_{\Gamma_{\mathrm{R},h}} \left( 2\dot{u}_h u_h + u_h^2 \nabla_{\mathrm{T}} \cdot V_h \right) dS \\ &= 2\int_{\Omega_h} \left( u_h (V_h \cdot \nabla) r + r\dot{u}_h + r u_h \nabla \cdot V_h \right) dV \\ &\quad - \int_{\Omega_h} \left( 2\nabla \dot{u}_h \cdot \nabla u_h - 2\nabla u_h \cdot (\nabla V_h) \nabla u_h + |\nabla u_h|^2 \nabla \cdot V_h \right) dV, \\ &\quad - \alpha \int_{\Gamma_{\mathrm{R},h}} \left( 2\dot{u}_h u_h + u_h^2 \nabla_{\mathrm{T}} \cdot V_h \right) dS, \end{aligned} \quad (4.23)$$

where expressions (4.18) and (4.20) have been used in the second equality. We recall that $\dot{u}_h \in W_{h,p}$. Thus, choosing $v_h = \dot{u}_h$ in equation (4.13) evaluated at $t = 0$, the three terms in expression (4.23) containing $\dot{u}_h$ cancel, which, together with the product rule, leaves us with

$$\begin{aligned} dJ_h(\Omega_h; V_h) &= 2\int_{\Omega_h} u_h \nabla \cdot (V_h r) \, dV + \int_{\Omega_h} \left( 2\nabla u_h \cdot (\nabla V_h) \nabla u_h - \nabla \cdot V_h |\nabla u_h|^2 \right) dV \\ &\quad - \alpha \int_{\Gamma_{\mathrm{R},h}} u_h^2 \nabla_{\mathrm{T}} \cdot V_h \, dS \\ &= 2\int_{\Gamma_{\mathrm{R},h}} \mathbf{n} \cdot V_h u_h r \, dS - \alpha \int_{\Gamma_{\mathrm{R},h}} u_h^2 \nabla_{\mathrm{T}} \cdot V_h \, dS \\ &\quad - \int_{\Omega_h} \left( 2 r V_h \cdot \nabla u_h - 2\nabla u_h \cdot (\nabla V_h) \nabla u_h + \nabla \cdot V_h |\nabla u_h|^2 \right) dV, \end{aligned} \quad (4.24)$$

where integration by parts and the fact that $u_h$ vanishes on $\Gamma_{\mathrm{D},h}$ has been used in the second equality.

The somewhat complicated expression (4.24) is called the *volume form* or the *weak form* of the directional derivative. This form is the one usually employed in implementations that aim for an exact derivative of the discrete objective function, and we note that the directional derivative depends in a quite intricate way on the transformation and the solution in the interior of $\Omega$. As shown in the previous publication [3], an analysis similar to the one done here holds also in the case before discretization. To achieve the simpler boundary form, or the *strong form* (4.12) of the directional derivative, integration by parts and some additional calculations have to be carried out, as exemplified in the previous publication [3]. We will perform these calculations also here in the discrete case, but the final expression will not be as simple as formula (4.12).

Let $\mathcal{S}_h^0 \subset \mathcal{S}_h$ be the set of element faces strictly inside the domain $\Omega$. For each $S \in \mathcal{S}_h^0$, there thus are two element $K_1, K_2 \in \mathcal{T}_h$ such that $\bar{S} = \bar{K}_1 \cap \bar{K}_2$. For $q_1 \in C^0(\bar{K}_1)$, $q_2 \in C^0(\bar{K}_2)$ and employing a notation of the type used for Discontinuous Galerkin methods (cf. [1, p. 1756]), we define the jump quantity

$$[\![q\mathbf{n}]\!]\big|_S = q_1\big|_S \mathbf{n}_1 + q_2\big|_S \mathbf{n}_2, \quad (4.25)$$

where $\mathbf{n}_1$ and $\mathbf{n}_2 = -\mathbf{n}_1$ are the normals on $S$ pointing outward from $K_1$ and $K_2$, respectively. We recall an integration-by-parts formula for vector-valued continuous functions $\boldsymbol{\psi}$ and scalar functions $q$ with jump discontinuities across element surfaces.[2] For $\boldsymbol{\psi} \in H^1(\Omega_h)^d$ and $q \in C^1(\bar{\mathcal{T}}_h)$ — that is, $q|_K \in C^1(\bar{K})$ for each $K \in \mathcal{T}_h$ — hold the integration-by-parts formula

$$\int_{\Omega_h} \nabla \cdot \boldsymbol{\psi} \, q \, dV + \sum_{K \in \mathcal{T}_h} \int_K \boldsymbol{\psi} \cdot \nabla q \, dV = \int_{\partial \Omega_h} \mathbf{n} \cdot \boldsymbol{\psi} \, q \, dS + \sum_{S \in \mathcal{S}_h^0} \int_S \boldsymbol{\psi} \cdot [\![\mathbf{n} q]\!] \, dS. \quad (4.26)$$

---

[2]As with notation (4.25), such formulas are standard tools in the context of Discontinuous Galerkin methods [1, equation (3.6)]

(A nonzero last term is a consequence of a presence of jump discontinuities in $q$.) Using formula (4.26) with $\boldsymbol{\psi} = \boldsymbol{V}_h$ and $q = |\nabla u_h|^2$, we find that

$$\int_{\Omega_h} \nabla \cdot \boldsymbol{V}_h |\nabla u_h|^2 \, dV = -\sum_{K \in \mathcal{T}_h} \int_K (\boldsymbol{V}_h \cdot \nabla) |\nabla u_h|^2 \, dV \\ + \int_{\partial \Omega_h} \boldsymbol{n} \cdot \boldsymbol{V}_h |\nabla u_h|^2 \, dS + \sum_{S \in \mathcal{S}_h^0} \int_S \boldsymbol{V}_h \cdot [\![\boldsymbol{n} |\nabla u_h|^2]\!] \, dS. \quad (4.27)$$

Substituting integration-by-parts formula (4.27) into expression (4.24), and recalling that $\boldsymbol{V}_h$ vanishes on $\Gamma_{\mathrm{D},h}$, we obtain

$$\begin{aligned} dJ_h(\Omega_h; \boldsymbol{V}_h) &= \int_{\Gamma_{\mathrm{R},h}} \boldsymbol{n} \cdot \boldsymbol{V}_h (2 u_h r - |\nabla u_h|^2) \, dS - 2 \int_{\Omega_h} r \boldsymbol{V}_h \cdot \nabla u_h \, dV \\ &\quad + \sum_{K \in \mathcal{T}_h} \int_K \underbrace{\left( (\boldsymbol{V}_h \cdot \nabla) |\nabla u_h|^2 + 2 \nabla u_h \cdot (\nabla \boldsymbol{V}_h) \nabla u_h \right)}_{= 2 \nabla u_h \cdot \nabla \left( (\boldsymbol{V}_h \cdot \nabla) u_h \right)} \, dV \\ &\quad - \sum_{S \in \mathcal{S}_h^0} \int_S \boldsymbol{V}_h \cdot [\![\boldsymbol{n} |\nabla u_h|^2]\!] \, dS - \alpha \int_{\Gamma_{\mathrm{R},h}} u_h^2 \nabla_{\mathrm{T}} \cdot \boldsymbol{V}_h \, dS \\ &= \int_{\Gamma_{\mathrm{R},h}} \boldsymbol{n} \cdot \boldsymbol{V}_h (2 u_h r - |\nabla u_h|^2) \, dS - \sum_{S \in \mathcal{S}_h^0} \int_S \boldsymbol{V}_h \cdot [\![\boldsymbol{n} |\nabla u_h|^2]\!] \, dS \\ &\quad + 2 \sum_{K \in \mathcal{T}_h} \int_K \left( \nabla u_h \cdot \nabla (\boldsymbol{V}_h \cdot \nabla) u_h - r \boldsymbol{V}_h \cdot \nabla u_h \right) \, dV - \alpha \int_{\Gamma_{\mathrm{R},h}} u_h^2 \nabla_{\mathrm{T}} \cdot \boldsymbol{V}_h \, dS. \end{aligned} \quad (4.28)$$

To continue, we need to apply integration by parts on the last term in expression (4.28). For this purpose, we will employ the following tangential divergence theorem for a function $\boldsymbol{\psi} \in C^1(\Gamma)^d$ defined on bounded smooth surface $\Gamma$ with a piecewise-smooth boundary $\partial \Gamma$,

$$\int_\Gamma \nabla_{\mathrm{T}} \cdot \boldsymbol{\psi} \, dS = \int_\Gamma \kappa \boldsymbol{n} \cdot \boldsymbol{\psi} \, dS + \int_{\partial \Gamma} \boldsymbol{m} \cdot \boldsymbol{\psi} \, dl, \quad (4.29)$$

where $\kappa = \nabla_{\mathrm{T}} \cdot \boldsymbol{n}$ is the summed curvature and $\boldsymbol{m}$ the outward-directed co-normal field on $\partial \Gamma$. Formula (4.29) follows from the standard vector-calculus Stokes theorem equaling the flux of the curl of a vector field $\boldsymbol{F}$ over $\Gamma$ to its circulation along $\partial \Gamma$. First note that the directed normal field $\boldsymbol{n}$ on $\Gamma$ can be extended into a tubular neighborhood of $\Gamma$ with the help of the signed distance function[3] $b_\Gamma$ chosen such that $\boldsymbol{n} = \nabla b_\Gamma|_\Gamma$. Formula (4.29) then follows from the choice $\boldsymbol{F} = \nabla b_\Gamma \times \boldsymbol{\psi}$.

Now let $S \in \mathcal{S}_h$ be any mesh face on $\Gamma_{\mathrm{R},h}$. Applying formula (4.29) with $\Gamma = S$ and $\boldsymbol{\psi} = -\boldsymbol{V}_h u_h^2$, together the product rule, we find that

$$-\int_S u_h^2 \nabla_{\mathrm{T}} \cdot \boldsymbol{V}_h \, dS = \int_S \boldsymbol{V}_h \cdot \nabla_{\mathrm{T}} u_h^2 \, dS - \int_S \kappa_S \boldsymbol{n} \cdot \boldsymbol{V}_h u_h^2 \, dS - \int_{\partial S} \boldsymbol{m} \cdot \boldsymbol{V}_h u_h^2 \, dl, \quad (4.30)$$

where we have introduced the *tangential gradient*, defined by

$$\nabla_{\mathrm{T}} f = \nabla f - \boldsymbol{n} \frac{\partial f}{\partial n}. \quad (4.31)$$

Since $S$ is a planar surface, $\kappa_S = 0$; however, we will keep the term in order to expose the parallel to formula (4.12). Moreover, keeping $\kappa_s$ will make the final formula correct also when using isoparametric elements with curved faces at the boundary. Summing equality (4.30) over all faces on $\Gamma_{\mathrm{R},h}$, we find the integration-by-parts formula

$$-\int_{\Gamma_{\mathrm{R},h}} u_h^2 \nabla_{\mathrm{T}} \cdot \boldsymbol{V}_h \, dS = \int_{\Gamma_{\mathrm{R},h}} \boldsymbol{V}_h \cdot \nabla_{\mathrm{T}} u_h^2 \, dS - \int_{\Gamma_{\mathrm{R},h}} \kappa_h \boldsymbol{n} \cdot \boldsymbol{V}_h u_h^2 \, dS - \sum_{e \in \mathcal{E}_h} \int_{e \cap \Gamma_{\mathrm{R},h}} [\![\boldsymbol{m}]\!] \cdot \boldsymbol{V}_h u_h^2 \, dl, \quad (4.32)$$

---

[3] Distance functions, signed distance functions, and tubular neighborhoods will be covered in more detail in § 5.1

where $\kappa_h$ almost everywhere on $\Gamma_{R,h}$ is defined by $\kappa_h|_S = \kappa_S$ for $S \in \mathscr{S}_h \cap \Gamma_{R,h}$, and where

$$[\![\boldsymbol{m}]\!] = \boldsymbol{m}_1 + \boldsymbol{m}_2 \tag{4.33}$$

in which $\boldsymbol{m}_1$ and $\boldsymbol{m}_2$ are the outward co-normals associated with the two surfaces $S_1$ and $S_2$ that meet at the line segent $e$; see the right picture in figure 1. By definition (4.31),

$$\begin{aligned}
\boldsymbol{V}_h \cdot \nabla_{\mathrm{T}} u_h^2 &= \boldsymbol{V}_h \cdot \nabla u_h^2 - \boldsymbol{n} \cdot \boldsymbol{V}_h \frac{\partial}{\partial n} u_h^2 \\
&= 2 u_h \boldsymbol{V}_h \cdot \nabla u_h - \boldsymbol{n} \cdot \boldsymbol{V}_h \frac{\partial}{\partial n} u_h^2,
\end{aligned} \tag{4.34}$$

which substituted into expression (4.32) yields

$$-\int_{\Gamma_{R,h}} u_h^2 \nabla_{\mathrm{T}} \cdot \boldsymbol{V}_h \, dS = 2 \int_{\Gamma_{R,h}} u_h \boldsymbol{V}_h \cdot \nabla u_h \, dS - \int_{\Gamma_{R,h}} \boldsymbol{n} \cdot \boldsymbol{V}_h \frac{\partial}{\partial n} u_h^2 \, dS \\ - \int_{\Gamma_{R,h}} \kappa_h \boldsymbol{n} \cdot \boldsymbol{V}_h u_h^2 \, dS - \sum_{e \in \mathscr{E}_h} \int_{e \cap \Gamma_{R,h}} [\![\boldsymbol{m}]\!] \cdot \boldsymbol{V}_h u_h^2 \, dl. \tag{4.35}$$

Substituting formula (4.35) into expression (4.28), we finally arrive at

$$\begin{aligned}
dJ_h(\Omega_h; \boldsymbol{V}_h) &= \\
&= \int_{\Gamma_{R,h}} \boldsymbol{n} \cdot \boldsymbol{V}_h \Big( 2r u_h - |\nabla u_h|^2 - \alpha \Big( \frac{\partial}{\partial n} u_h^2 + \kappa_h u_h^2 \Big) \Big) dS - \sum_{e \in \mathscr{E}_h} \alpha \int_{e \cap \Gamma_{R,h}} [\![\boldsymbol{m}]\!] \cdot \boldsymbol{V}_h u_h^2 \, dl \\
&\quad - \sum_{S \in \mathscr{S}_h^0} \int_S \boldsymbol{V}_h \cdot [\![\boldsymbol{n} |\nabla u_h|^2]\!] \, dS \\
&\quad + 2 \sum_{K \in \mathscr{T}_h} \int_K \Big( \nabla u_h \cdot \nabla (\boldsymbol{V}_h \cdot \nabla) u_h - r \boldsymbol{V}_h \cdot \nabla u_h \Big) dV + 2\alpha \int_{\Gamma_{R,h}} u_h (\boldsymbol{V}_h \cdot \nabla) u_h \, dS.
\end{aligned} \tag{4.36}$$

Note that this directional derivative expressions is completely equivalent to the "volume form" (4.24), although the relation to the before-discretization boundary expression (4.12) is clearer here.

## 5 Domain paths using surface dilations

### 5.1 Delfour dilations

In 2018, Michel Delfour [10] introduced a way to construct domain paths $t \mapsto \Omega_t$ that is fundamentally different from the domain transformation approach considered so far. Instead of transformations, the domain paths are built from *dilations* of a lower-dimensional object $E$ into $\mathbb{R}^d$. Delfour uses such dilations to generalize the concept of *topological derivative*, that is, when sensitivities are calculated with respect to the introduction of a small hole in the domain.

Thus, let $E$ be a closed subset of $\mathbb{R}^d$ of dimension $m = 0, \ldots, d-1$ with finite $m$-dimensional measure. For instance, $E$ could be a point ($m = 0$), a curve segment ($m = 1$), or a surface patch ($m = 2$). The $r$-*dilated set* $E_r$ is the set of points $\boldsymbol{x} \in \mathbb{R}^d$ within a distance $r > 0$ from $E$, that is, the set of points $\boldsymbol{x}$ such that $d_E(\boldsymbol{x}) = \inf_{\boldsymbol{e} \in E} |\boldsymbol{x} - \boldsymbol{e}| \leq r$. (Note that $E$ and $E_r$ are of different dimensions; $E$ is of dimension $m < d$, whereas the dilated $E_r$ is of dimension $d$). The dilations $E_r$ will be be used to construct domain paths $\Omega_t$ with parameter $t$ being the volume of the ball of codimension $m$ with radius $r$. That is, $t = \alpha_{d-m} r^{d-m}$, where $\alpha_{d-m}$ is the volume of the unit ball of codimension $m$. Various perturbations of $\Omega$ may be defined using $E_r$. Here we consider the domain path $t \mapsto \Omega_t = \Omega \setminus E_{r(t)}$, that is, the domain $\Omega$ with the dilated region cut out. The idea is now to define directional semiderivatives of a shape function $J(\Omega)$ with respect to $E$ as the limit of the difference quotient

$$dJ(\Omega; E) = \lim_{t \to 0^+} \frac{J(\Omega_t) - J(\Omega)}{t}, \tag{5.1}$$

if the limit exists. If $E \subset \Omega$ and $J(\Omega) = \int_\Omega dV$, the quantity $-dJ(\Omega; E)$ turns out to be the *Minkowski content* of $E$, which equals the Hausdorff measure[4] of $E$ for so-called *m-rectifiable* and closed $E$.

As a first example, let $m = 0$ and $E = \{x_0\}$ for some $x_0 \in \Omega$ [10, Ex. 4.1]. The dilated region $E_r$ is a ball of radius $r$ centered at $x_0$, $t = \alpha_d r^d = |E_r|$, and $\Omega_t = \Omega \setminus E_r$. For $J_1$ from expression (3.1), it then follows that

$$dJ_1(\Omega; \{x_0\}) = \lim_{t \to 0^+} \frac{J(\Omega_t) - J(\Omega)}{t} = \lim_{r \to 0^+} \frac{1}{|E_r|} \left( \int_{\Omega \setminus E_r} f \, dV - \int_\Omega f \, dV \right)$$
$$= -\lim_{r \to 0^+} \frac{1}{|E_r|} \int_{E_r} f \, dV = -f(x_0), \tag{5.2}$$

which equals the usual topological derivative of $J_1$ for a fixed integrand $f \in C(D)$.

Similarly, for $m = 1$, if $E = \gamma$, where $\gamma$ is a closed $C^2$ curve embedded in $\Omega$, an analogous calculation [10, Ex. 4.16] shows that, for $f \in W^{1,1}(D)$,

$$dJ_1(\Omega; \gamma) = -\int_\gamma f \, ds. \tag{5.3}$$

However, the focus here will be the case $m = d - 1$ and $E = \partial \Omega$, which, as we will see, leads to a shape derivative. According to the above scheme, we should have $t = \alpha_{d-m} r^{d-m} = 2r$. However, $E$ is in this case not a subset of $\Omega$, as assumed above. Indeed, $E_r = U_r(\partial \Omega)$, where

$$U_r(\partial \Omega) = \left\{ x \in \mathbb{R}^d \mid d_{\partial \Omega}(x) \leq r \right\} \tag{5.4}$$

is the $r$-tubular neighborhood of $\partial \Omega$, which consists of two parts, one inside and one outside $\Omega$. It will only be the interior part $E_r \cap \Omega$ that will perturb $\Omega$ when forming $\Omega \setminus E_r$. Recall the property that for rectifiable $E \subset \Omega$ and $J(\Omega) = \int_\Omega dV$, the quantity $-dJ(\Omega; E)$ is the measure of $E$. In order to retain this property when $E = \partial \Omega$, we need to choose $t = r$, otherwise $-dJ(\Omega; \partial \Omega)$ will only represent half of the measure of $\partial \Omega$.

Now, forming the difference quotient for shape function $J_1$ with $f \in W^{1,1}(D)$, we obtain, since $\Omega_t = \Omega \setminus U_t$,

$$dJ_1(\Omega; \partial \Omega) = \lim_{t \to 0^+} \frac{1}{t} \left( \int_{\Omega_t} f \, dV - \int_\Omega f \, dV \right) = -\lim_{t \to 0^+} \frac{1}{t} \int_{U_t \cap \Omega} f \, dV = -\int_{\partial \Omega} f \, dS. \tag{5.5}$$

Thus, the directional derivative of $J_1$ with respect to dilation of the domain boundary coincides with the directional derivative (4.7a) with velocity field $V = -n$. (Integrand $f$ is fixed, so $f' = 0$.)

Moreover, proceeding similarly as Delfour [10, § 5], for the choice $E = \partial \Omega$ and assuming sufficient smoothness, it follows that

$$dJ_2(\Omega; \partial \Omega) = -\int_{\partial \Omega} \left( \frac{\partial f}{\partial n} + \kappa f \right) dS. \tag{5.6}$$

Again, the directional derivative of $J_2$ with respect to dilation of the domain boundary coincides with the directional derivative (4.7b) with the velocity field $V = -n$. For shape optimization applications, the directional derivatives (5.5) and (5.6) are of limited use, since only uniform extensions and contraction of the domain are addressed. To achieve a localized control over boundary shapes, we will reformulate and generalize the dilation concept using level-set functions.

## 5.2 Dilation through level sets

In the case $m = d - 1$ and $E = \partial \Omega$, the Delfour dilation generated the domain path $t \to \Omega_t = \Omega \setminus U_t$, where $U_t$ is the tubular neighborhood of $\partial \Omega$, defined in expression (5.4). Equivalently, this domain path can be constructed through the use of the *signed distance function* $b_\Omega : D \to \mathbb{R}$, defined as

$$b_\Omega(x) = d_\Omega(x) - d_{D \setminus \Omega}(x), \tag{5.7}$$

---

[4] scaled such that that it agrees with the $m$-dimensional Lebesgue measure

where, as before, $d_E(\boldsymbol{x}) = \inf_{\boldsymbol{e} \in E} |\boldsymbol{x} - \boldsymbol{e}|$. The signed distance function is negative and positive inside and outside of $\Omega$, respectively, and its absolute value specifies the distance to the boundary $\partial \Omega$. In terms of the signed distance function, $\Omega_t = \Omega \setminus U_t$ can alternatively be defined as

$$\Omega_t = \{\, \boldsymbol{x} \in D \mid b_\Omega(\boldsymbol{x}) + t \leq 0 \,\}. \tag{5.8}$$

Aiming to generalize and localize characterization (5.8), let $\phi : D \to \mathbb{R}$ be a level-set function, that is, a Lipschitz-continuous function $\phi$ that partitions $D$ such that

$$\Omega = \{\, \boldsymbol{x} \in D \mid \phi(\boldsymbol{x}) < 0 \,\}, \quad \partial\Omega = \{\, \boldsymbol{x} \in D \mid \phi(\boldsymbol{x}) = 0 \,\},$$
$$D \setminus \overline{\Omega} = \{\, \boldsymbol{x} \in D \mid \phi(\boldsymbol{x}) > 0 \,\}. \tag{5.9}$$

The signed distance function is an example of a level-set function. Moreover, let $w : D \to [0, 1]$ be a Lipschitz function with local support $\operatorname{supp} w \subset\subset D$, such that $w > 0$ in the interior of its support. The function $w$ will later be chosen as a continuous, piecewise-linear finite-element basis function. Definition (5.8) can then be generalized to

$$\Omega_t(w) = \{\, \boldsymbol{x} \in D \mid \phi(\boldsymbol{x}) + tw(\boldsymbol{x}) \leq 0 \,\}. \tag{5.10}$$

We will define, in the discretized case, semi-derivatives of shape functions (3.1) using domain paths as in definition (5.10). Note that this definition does not involve any domain transformation!

### 5.3 Fictitious-domain shape optimization and dilations

Instead of requiring the computational mesh to conform to changing domain shapes, an alternative is offered by *fictitious-domain*, also called *domain embedding*, methods. There are many flavors of this idea, but common to all of them is that a fixed computational mesh is introduced on the hold-all $D$, and that the domain boundary $\partial\Omega$ cuts through $D$. In the most general versions of the method — the kind considered here — the boundary is also allowed to intersect the interior of the mesh cells. In the context of shape optimization, the big advantage with this class of methods is that mesh deformation strategies are not needed. The cost of the method is that additional efforts are usually needed to assure numerical stability and accurate imposition of boundary conditions.

As opposed to the body-fitted methods discussed previously, the fictitious-domain approach does not fit well with the domain-transformation concept; there are no domain transformations naturally involved when moving boundary $\partial\Omega$ over a fixed background mesh. However, the dilation technique introduced in § 5.2 is a good framework for shape calculus in the context of fictitious-domain methods.

Using a similar notation as in § 4.2, we assume a simplicial triangulation $\mathcal{T}_h$, now of the full hold-all $D$, and a function space $W_{h,1}$ of continuous, piecewise-linear functions on this triangulation. The boundary of each admissible domain will be defined as the zero level set of a function $\phi \in W_{h,1}$, and perturbations of the domain are defined through a perturbation of the level-let function,

$$\phi_t(\boldsymbol{x}) = \phi(\boldsymbol{x}) + tw(\boldsymbol{x}), \tag{5.11}$$

where $t \geq 0$ and $w$ a Lagrangian basis function for $W_{h,1}$. We are thus in the framework of § 5.2.

We start with the directional semiderivative of the volume integral $J_1$ in expression (3.1). The domain transformation approach used to show expression (4.7a) required $f \in W^{1,1}(D)$. For the proof of the corresponding formula in the current framework, we can considerably relax this condition. We only need to assume that the integrand is given in terms of functions $f(t) \in C^0(\overline{\mathcal{T}}_h)$, that is, $f(t)|_K \in C^0(\overline{K})$ for each $K \in \mathcal{T}_h$. Thus, jump discontinuities between elements of the mesh are permitted. In this case, the following theorem holds.

**Theorem 5.1.** *Under perturbation* (5.11) *and for* $t \mapsto f(t)$ *and* $t \mapsto f'(t)$ *continuous in some nonempty interval* $[0, t_{max}]$ *such that* $f(t), f'(t) \in C^0(\overline{\mathcal{T}}_h)$ *on* $(0, t_{max})$, *the directional semiderivative of volume integral*

$$J_1(\phi_t) = \int_{\Omega_t} f_t \, dV \tag{5.12}$$

*at $t = 0$ satisfies*

$$dJ_1(\phi; w) = \lim_{t \to 0^+} \frac{1}{t} \big(J_1(\phi_t) - J_1(\phi)\big) = \int_\Omega f' \, dV - \int_{\partial\Omega} f \frac{w}{|\partial_n \phi|} \, dS. \tag{5.13}$$

Theorem 5.1 is an immediate consequence of theorem 6.7 in § 6, which proves the statement for a fixed $f \in C^0(\bar{\mathscr{T}}_h)$. When $\partial\Omega$ only cuts through the interior of the mesh elements, as in the right picture of figure 2, there is no ambiguity in formula (5.13). In fact, the limits of $t \to 0^+$ and $t \to 0^-$ will then agree. However for a domain like in the left picture of figure 2, when $\int_{S \cap \partial\Omega} dS > 0$ for some mesh face $S$, both $f$ and $\partial_n \phi$ are typically discontinuous across $S$. Formula (5.13) then holds for the limit of these values from the *interior* of $\Omega$. Note that this ambiguity means that whenever $\int_{S \cap \partial\Omega} dS > 0$ for some mesh face $S$, $J_1$ possesses only a one-sided derivative; when the limit $t \to 0^-$ is considered, it will be the limit of the values from the exterior that should be employed.

Regarding the shape calculus associated with surface integral $J_2$ of expression (3.1), note that any nonempty $K \cap \partial\Omega$ is a line segment for $d = 2$ and either a triangle or a quadrilateral for $d = 3$ (figure 4). For any such $K \cap \partial\Omega$, we will need the concept of co-normals, already introduced in § 4.2, illustrated for the current case in figure 6. These will be used below to define limits and jumps across $\partial\Omega \cap S$ where $S$ is a mesh surface. Such a $\partial\Omega \cap S$ is a line segment for $d = 3$ and a point for $d = 2$.

The following theorem is an immediate consequence of theorem 6.11 in § 6, which proves the statement for a fixed $f \in \mathscr{C}^1(\bar{\mathscr{T}}_h)$. As for theorem 5.1, nondifferentiability can be expected when the boundary $\partial\Omega$ is partly aligned with a mesh face, as in the left picture of figure 2. In fact, the nondifferentiability situation is even more severe here due to the last term in the expression for $dJ_2$ below. Possible ambiguities with respect to this term is not as easy to resolve as for $dJ_1$, which is why the proof of theorem 6.11 is restricted to domains where the boundary does not intersect any mesh nodes, such as the one exemplified to the right in figure 2.

**Theorem 5.2.** *Assume that the boundary $\partial\Omega$ of the domain $\Omega$ does not intersect with any mesh nodes of $\mathscr{T}_h$. Then, under perturbation (5.11) and for $t \mapsto f(t)$ and $t \mapsto f'(t)$ being continuous in some nonempty interval $[0, t_{max}]$ such that $f(t), f'(t) \in C^1(\bar{\mathscr{T}}_h)$ on $(0, t_{max})$, the directional derivative of surface integral*

$$J_2(\phi_t) = \int_{\partial\Omega_t} f_t \, dS \tag{5.14}$$

*at $t = 0$ satisfies*

$$dJ_2(\phi, w) = \int_{\partial\Omega} \left( f' - \frac{\partial f}{\partial n} \frac{w}{|\partial_n \phi|} \right) dS - \sum_{S \in \mathscr{S}_h} \int_{\partial\Omega \cap S} \boldsymbol{n}^S \cdot [\![ f \boldsymbol{m} ]\!] \frac{w}{|\partial_{n^S} \phi|} \, d\gamma, \tag{5.15}$$

*where $\boldsymbol{n}^S$ is the normal vector to $\partial\Omega \cap S$, located in $S$ and outward-directed from $\Omega$, and where $[\![ f \boldsymbol{m} ]\!] = f_1 \boldsymbol{m}_1 + f_2 \boldsymbol{m}_2$. Here, for $k = 1, 2$, $\boldsymbol{m}_k$ are the co-normals to $\partial\Omega \cap K_k$ at $\partial\Omega \cap S$, where $K_1, K_2 \in \mathscr{T}_h$ such that $\bar{S} = \bar{K}_1 \cap \bar{K}_2$. Each $\boldsymbol{m}_k$ lies in the plane of and is directed outward from $\partial\Omega \cap K_k$, and $f_k$ is the limit $f$ on $S$ defined by $f_k(\boldsymbol{x}) = \lim_{\epsilon \to 0^+} f(\boldsymbol{x} - \epsilon \boldsymbol{m}_k)$, for any $\boldsymbol{x} \in S$.*

*Remark 5.3.* For $d = 2$, $\partial\Omega \cap S$ is a point, so the the second integral in expression (5.15) should be interpreted as a point evaluation at $\partial\Omega \cap S$; that is, the expression simply is

$$dJ_2(\phi, w) = \int_{\partial\Omega} \left( f' - \frac{\partial f}{\partial n} \frac{w}{|\partial_n \phi|} \right) dS - \sum_{S \in \mathscr{S}_h} \boldsymbol{n}^S \cdot \left( [\![ f \boldsymbol{m} ]\!] \frac{w}{|\partial_{n^S} \phi|} \right) \bigg|_{\partial\Omega \cap S}, \tag{5.16}$$

and $\boldsymbol{n}^S$ is a unit vector along the mesh edge $S$. □

We note that expressions (5.13) and (5.15) constitute discrete versions of the directional derivative expressions (4.7). The weight $w/|\partial_n \phi|$ corresponds to the velocity $-\boldsymbol{n} \cdot \boldsymbol{V}$ and the last term in formula (5.15) corresponds to the term involving the boundary curvature in formula (4.7b). As will be seen from the proof of theorem 6.11, the first integral of expression (5.15) also formally includes a term $\kappa f$, as in formula (4.7b), which however vanishes since, in this case, $\partial\Omega \cap K$ is planar for each $K \in \mathscr{T}_h$.

### 5.4 Shape calculus of the model problem for fictitious-domain methods

On the triangulation $\mathscr{T}_h$ of the full hold-all $D$, we consider the space $W_{h,p}$ of continuous finite-element functions that are polynomials of maximal degree $p \geq 1$ on each $K \in \mathscr{T}_h$. Any such function can be expressed as a sum over Lagrangian basis functions $N_i^p$,

$$u_h(\boldsymbol{x}) = \sum_{i=1}^{I} u_i N_i^p(\boldsymbol{x}), \tag{5.17}$$

where $I$ is the number of evaluation nodes for the basis. (Function $N_i^1$ corresponds to the linear basis function $w$ used for the level-set functions.)

Now assume that the connected reference domain $\Omega_h \subset\subset D$ is defined as the set of points where level-set finite-element function $\phi \in W_{h,1}$ is strictly negative. Moreover, we consider perturbation (5.11), which generates a domain path $t \mapsto \Omega_{h,t}$ that cuts through the background mesh $\mathscr{T}_h$.

Analogously with the treatments in § 4.1 and 4.2, we will keep the Dirichlet boundary $\Gamma_D$ of problem (2.1) fixed under domain perturbations. Therefore, we adapt the mesh so that $\Gamma_D$ aligns with faces of the elements in $\mathscr{T}_h$. Then the homogeneous boundary condition $\Gamma_D$ is easy to assign just by removing the evaluation nodes on $\overline{\Gamma}_D$ from the expansion (5.17). To keep $\Gamma_D$ fixed, we only permit perturbations (5.11) for basis functions $w$ whose support does not intersect with $\Gamma_D$.

*Remark* 5.4. If $\Gamma_D$ is not aligned with the faces of $\mathscr{T}_h$, a more elaborate strategy is needed to assign the Dirichlet boundary condition, for instance using a Nitsche-type approach [8].

A fictitious-domain finite-element approximation of variational problem (4.8) can then be formulated as

$$u_h(t) \in \widehat{W}_{h,p}(t) \text{ such that}$$
$$\int_{\Omega_{h,t}} \nabla v_h \cdot \nabla u_h(t)\, dV + \alpha \int_{\Gamma_{R,h,t}} v_h u_h(t)\, dS = \int_{\Omega_{h,t}} r v_h \, dV \qquad \forall v_h \in \widehat{W}_{h,p}(t), \quad (5.18)$$

where here $\widehat{W}_{h,p}(t)$ is the space of restriction on $\Omega_{h,t}$ of functions in $W_{h,p}$. The discrete objective function is

$$J_h(\phi) = \int_{\Omega_{h,t}} r u_h(t)\, dV. \quad (5.19)$$

The approximation properties of problem (5.18) is similar to the mesh conforming discretization (4.13), but the condition number of the matrix with elements $\int_{\Omega_{h,t}} \nabla N_i^p \cdot \nabla N_j^p \, dV$ can be arbitrarily high for cases when $\Omega_{h,t} \cap K$ has a small measure for some $K \in \mathscr{T}_h$. Approximation (5.18) can be viewed as being in the class of cut-finite-element methods [7], for which a condition-number stabilization scheme, suggested by Burman [6], has been developed. In this scheme, an additional so-called ghost-penalty term is added to the variational form to bound some of the jumps in the normal derivatives across neighboring elements. As analyzed by Bernland et al. [4], the presence of the ghost penalty term affects the shape calculus only in exceptional cases, which is why we for simplicity ignore stabilization issues here.

In § 4.2, we concluded that it is the *material* derivative (4.4) of the finite-element function that stays in the finite-element space when mesh deformations are employed, whereas the shape derivative will not. The reason is that the support of each basis function is also transformed in that case. However, the mesh is fixed in the fictitious-domain case, so it is only the coefficients of expansion (5.17) that can change with domain-path parameter $t$. Thus, here it is the *shape* derivative that will stay in the finite-element space, whereas the material derivative is not even well defined, since no transformations are involved a priori.

Thus, the shape calculus of our model problem under fictitious-domain discretization can be carried out in the same manner as in § 4.1, but using formulas (5.13) and (5.15) instead of formulas (4.7), resulting in

$$dJ_h(\phi; w) = -\int_{\Gamma_{R,h}} \frac{w}{|\partial_n \phi|}\left(2r u_h - |\nabla u_h|^2 - \alpha \frac{\partial}{\partial n} u_h^2\right) dS + \alpha \sum_{S \in \mathscr{S}_h} \int_{\Gamma_{R,h} \cap S} \boldsymbol{n}^S \cdot [\![\boldsymbol{m}]\!] u_h^2 \frac{w}{|\partial_n \phi|}\, d\gamma. \quad (5.20)$$

We note that directional derivative (5.20) corresponds to the boundary expression (4.12) with the substitution

$$\boldsymbol{n} \cdot V = -\frac{w}{|\partial_n \phi|} \quad (5.21)$$

and with the jump of the co-normals $\boldsymbol{m}$ serving as a discrete measure of the curvature $\kappa$ of $\Gamma_R$. As seen in the proof of theorem 6.11, the summed curvature of each $K \cap \partial \Omega_h$ also appear in the derivation but vanishes due to the piecewise planar shape of the cut boundary.

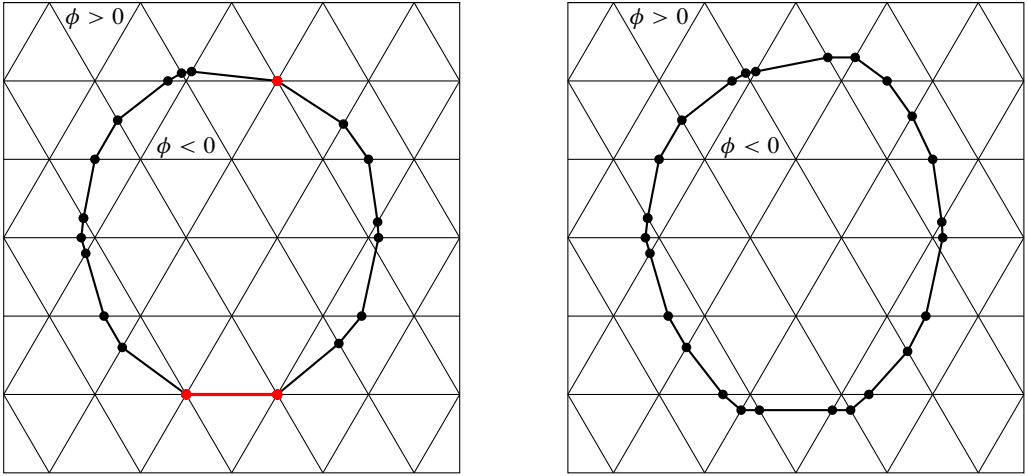

FIGURE 2. Examples for $d = 2$ of a rectangular hold-all $D$ and domains $\Omega$ generated by level-set functions $\phi$ defined on the mesh in the figure. Assumption 6.4 is satisfied for the domain to the right but not for the one to the left.

## 6 Proofs for the discrete fictitious-domain case

This section provides the precise assumptions and supplies the proofs behind theorems 5.1 and 5.2. For completeness and to make this section independently readable, there is some repetition of previously introduced notation and definitions.

Assume that the hold-all $D \subset \mathbb{R}^d$, for $d = 2$ or 3, is a simply connected bounded polyhedron ($d = 3$) or polygon ($d = 2$). We introduce a nondegenerate triangulation $\mathcal{T}_h$ of $D$ using simplices (tedrahedra for $d = 3$ and triangles for $d = 2$). Associated with $\mathcal{T}_h$, we consider the set of faces $\mathcal{S}_h$ and and subfaces $\mathcal{E}_h$. These are simplices of codimension 1 and 2, respectively. The faces $\mathcal{S}_h$ are mesh surfaces for $d = 3$ and mesh edges for $d = 2$, and the subfaces $\mathcal{E}_h$ are mesh edges for $d = 3$ and mesh points for $d = 2$. Associated with this mesh, we consider the space $W_h = W_{h,1}$ of continuous functions on $D$ that are linear on each $K \in \mathcal{T}_h$.

Now assume that $\phi \in W_h$ partitions $D$ such that

$$\Omega = \{ x \in D \mid \phi(x) < 0 \}, \quad \partial \Omega = \{ x \in D \mid \phi(x) = 0 \},$$
$$D \setminus \overline{\Omega} = \{ x \in D \mid \phi(x) > 0 \}, \tag{6.1}$$

and, moreover, that $\Omega$ is nonempty and connected. With this choice of $W_h$, the boundary $\partial \Omega$ comprises a set of polygons of codimension 1. Figure 2 illustrate example cases for $d = 2$.

*Remark* 6.1. Here we use the convention that each $K \in \mathcal{T}_h$, each $S \in \mathcal{S}_h$, and, for $d = 3$, each $E \in \mathcal{E}_h$ are relatively open in the sense that they exclude their boundaries $\partial K$, $\partial S$, and $\partial E$.

Now consider the perturbation

$$\phi_t(x) = \phi(x) + t w(x), \tag{6.2}$$

where $t \geq 0$ and $w$ is a Lagrangian nodal basis function for $V_h$, and the associated family of domains given by

$$\Omega_t = \{ x \in D \mid \phi_t(x) < 0 \}. \tag{6.3}$$

Linear Lagrangian basis functions satisfy $w \geq 0$, so $t'w(x) \leq tw(x)$ for $0 \leq t' \leq t$. Thus, if $x \in \Omega_t$, then $\phi(x) + t'w(x) \leq \phi(x) + tw(x) < 0$, so $x \in \Omega_{t'}$ also. Hence, $\Omega_t \subset \Omega_{t'}$, for $t' \leq t$; that is, the domains are shrinking for increasing positive $t$.

*Remark* 6.2. The choice $t \geq 0$ and the limits $t \to 0^+$ considered below is arbitrary. We could as well consider $t \leq 0$ and limits $t \to 0^-$. The domains would then be expanding for decreasing negative $t$, but the analysis below would be completely analogous also in this case. □

Since in the end, we will only be concerned with limits as $t \to 0^+$, we may without loss of generality restrict the perturbation to be small enough for the following condition to hold.

*Assumption* 6.3. There is a $t_{\max} > 0$ such that for each $K \in \mathscr{T}_h$, if $\partial \Omega_t \cap K$ is either empty or non empty for some $t \in (0, t_{\max}]$, it has the same property for each $t \in (0, t_{\max}]$. □

That is, it cannot happen that $\partial \Omega_t$ starts intersecting an element $K$ only for a large enough $t$ in $(0, t_{\max}]$. Assumption 6.3 yields some immediate consequences:

1. Since $\phi_t$ is a continuous function, assumption 6.3 will hold also for the elements in $\mathscr{S}_h$ and $\mathscr{E}_h$. That is, for each $S \in \mathscr{S}_h$ and for each $e \in \mathscr{E}_h$, if $\partial \Omega_t \cap S$ or $\partial \Omega_t \cap e$ is either empty or non empty for some $t \in (0, t_{\max}]$, the same property holds for each $t \in (0, t_{\max}]$.

2. No topology changes of $\Omega_t$ can happen as $t$ is varied, only boundary perturbations of the initial domain.

3. The center node $\boldsymbol{x}_w$ of basis function $w$ in perturbation (6.2), that is, the mesh vertex satisfying $w(\boldsymbol{x}_w) = 1$, will not intersect $\partial \Omega_t$ for any $t \in (0, t_{\max}]$, since $\boldsymbol{x}_w \notin K$ for any $K \in \mathscr{T}_h$ ($K$ is open; remark 6.1) and such an intersection would therefore contradict assumption 6.3.

Thus, item 3 says that $\phi_t(\boldsymbol{x}_w) \neq 0$ for $t \in (0, t_{\max}]$. Note, however, that it may happen that $\phi(\boldsymbol{x}_w) = 0$; that is, the boundary $\partial \Omega$ may intersect nodes in the mesh, as exemplified with the mesh nodes marked in red in figure 2.

For $k \in \mathbb{N}$, let $C^k(\overline{\mathscr{T}}_h)$ denote the space of functions such that $f|_K \in C^k(\overline{K})$ for each $K \in \mathscr{T}_h$. In practice, these will be finite-element functions $f$ that are polynomials on each element, but which are allowed to possess jump discontinuities between elements. For integrals

$$J_1(\phi) = \int_{\Omega(\phi)} f \, dV, \qquad J_2(\phi) = \int_{\partial \Omega(\phi)} f \, dS, \qquad (6.4)$$

containing such function, we will devise expressions for the directional semiderivatives

$$dJ_k(\phi; w) = \lim_{t \to 0^+} \frac{1}{t}\big(J_k(\phi_t) - J_k(\phi)\big), \quad k = 1, 2. \qquad (6.5)$$

Theorems 6.7 and 6.11 below are general enough to cover objective functions with such integrands $f$, a generality that will be needed when applying the formulas to finite-element variational expressions. However, there are issues with this generality for domains such as the one to the left in figure 2, where $\partial \Omega$ intersects nodes in the mesh and, in particular, partly aligns with the mesh, instances marked red in the figure. The seriousness of this issue differs for $dJ_1$ and $dJ_2$, being worse for the latter.

The directional semiderivate $dJ_1$ contains integrals over $\partial \Omega$ with objects such as $f$ and $\partial_n \phi$ possessing jump discontinuities at mesh surfaces. Limiting values of these are well defined almost everywhere on $\partial \Omega$ — but in general different — when approaching opposite sides of a mesh surface. A consequence is that the directional semiderivatives of $J_1$ as $t \to 0^+$ and as $t \to 0^-$ will differ in general, meaning that only semidifferentiability has any chance of hold for such domains.

Semiderivative $dJ_2$ includes in addition to integrals over $\partial \Omega$ also line integrals ($d = 3$) or point evaluations ($d = 2$) of quantities evaluated at points where the boundary intersects a mesh surface. These terms will be well defined whenever these intersections are strictly inside the mesh surface, but it is less clear how to treat the general case. Thus, we will introduce an additional requirement, assumed for $dJ_2$ but optional for $dJ_1$; when assumed, it will force the directional semiderivatives of $J_1$ as $t \to 0^+$ and as $t \to 0^-$ to agree.

*Assumption* 6.4. The boundary $\partial \Omega$ does not intersect any mesh points of the triangulation $\mathscr{T}_h$.

Figure 2 illustrate domains satisfying and not satisfying the assumption.

### 6.1 The domain integral $J_1$

For $t \in (0, t_{\max}]$, it holds for domain integral $J_1$ in definition (6.4) that

$$\frac{1}{t}\big(J_1(\phi_t) - J_1(\phi)\big) = \frac{1}{t} \int_{\Omega_t} f \, dV - \frac{1}{t} \int_{\Omega} f \, dV = -\frac{1}{t} \int_{E_t} f \, dV, \qquad (6.6)$$

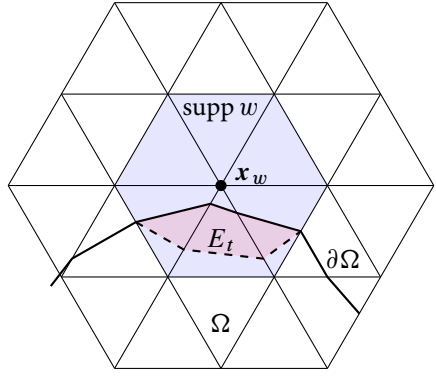

FIGURE 3. An example of a region $E_t = \Omega \setminus \overline{\Omega}_t$ using level-set perturbation $\phi_t = \phi + tw = 0$, where the support of the linear Lagrangian basis function w is indicated with a tinted blue color.

with $E_t = \Omega \setminus \overline{\Omega}_t$, which is nonempty and open provided that the support of basis function $w$ possesses a nonempty intersection with $\Omega$. Figure 3 shows an example for $d = 2$. We will parametrize region $E_t$, element by element, using functions $\widehat{X}$ of the following lemma. Roughly speaking, each nonempty $E_t \cap K$ will be parametrized along rays crossing the node $\boldsymbol{x}_w$ of basis function $w$ and the face $S$ of $K$ positioned opposite to $\boldsymbol{x}_w$.

**Lemma 6.5.** *Consider perturbation (6.2) for some $t \in (0, t_{max}]$, according to assumption 6.3. Let $K \in \mathcal{T}_h$ be an element in the support of basis function $w$. Then, for some $S \in \mathcal{S}_h$, there exists a line segment ($d = 2$) or a polygonal area ($d = 3$) $\widehat{S} \subset S$ and a diffeomorphism $\widehat{X} : \widehat{S} \times (0, t)$ such that $\widehat{X}(\widehat{S}, (0, t)) = E_t \cap K$ and $\widehat{X}(\widehat{S}, \tau) = \partial\Omega_\tau \cap K$ for each $\tau \in (0, t)$.*

*Proof.* We may assume that $E_t \cap K$ is nonempty, otherwise the statement is vacuously true. Denote by $\boldsymbol{x}_w$ the center node of basis function $w$, that is, the node satisfying $w(\boldsymbol{x}_w) = 1$. Note that $\boldsymbol{x}_w$ is one of the vertices of $K$. Let $S \in \mathcal{S}_h$ be the face of $K$ opposite to $\boldsymbol{x}_w$, that is, the unique face $S \subset \overline{K}$ such that $\boldsymbol{x}_w \notin \overline{S}$. Each point $\boldsymbol{x} \in K$ lies somewhere on a ray from $S$ through $\boldsymbol{x}_w$. More precisely, for each $\boldsymbol{x} \in K$, there is a unique pair $(\boldsymbol{x}_S, \sigma) \in S \times (0, 1)$ such that

$$\boldsymbol{x} = \boldsymbol{x}_S + \sigma(\boldsymbol{x}_w - \boldsymbol{x}_S). \tag{6.7}$$

In particular, since $\partial\Omega_\tau \cap K$ is nonempty for each $\tau \in (0, t)$ due to assumption 6.3, given any $\boldsymbol{x}_\tau \in \partial\Omega_\tau \cap K$, there is a unique pair $\boldsymbol{x}_S \in S$ and $\sigma \in (0, 1)$ such that

$$\boldsymbol{x}_\tau = \boldsymbol{x}_S + \sigma(\boldsymbol{x}_w - \boldsymbol{x}_S). \tag{6.8}$$

We will now specify the domain $\widehat{S} \subset S$ for $\boldsymbol{x}_S$ above. Figure 4 illustrates, for $d = 3$, the conceptual geometry options for $\widehat{S}$. Formula (6.8) specifies $\partial\Omega_\tau \cap K$ in terms of the zero set of $\phi_\tau|_K$, which means that $\boldsymbol{x}_S$ and $\boldsymbol{x}_w$ must be on different sides of $\partial\Omega_\tau \cap K$, and thus that $\phi_\tau(\boldsymbol{x}_S)$ and $\phi_\tau(\boldsymbol{x}_w)$ are of opposite signs, that is, $\phi_\tau(\boldsymbol{x}_S)\phi_\tau(\boldsymbol{x}_w) < 0$. This inequality suggests the following characterization of the set of points $\widehat{S} \subset S$ used to generate $\partial\Omega_\tau \cap K$ according to formula (6.8),

$$\widehat{S} = \{ \boldsymbol{x} \in S \mid \phi(\boldsymbol{x})\phi_{t_{\max}}(\boldsymbol{x}_w) < 0 \}, \tag{6.9}$$

where we have used that $\mathrm{sgn}(\phi_\tau(\boldsymbol{x}_S)\phi_\tau(\boldsymbol{x}_w)) = \mathrm{sgn}(\phi(\boldsymbol{x})\phi_{t_{\max}}(\boldsymbol{x}_w))$, since $\phi_\tau = \phi$ on $S$, due to the vanishing of $w$ on $S$, and since $\mathrm{sgn}\,\phi_\tau = \mathrm{sgn}\,\phi_{t_{\max}}$, due to assumption 6.3.

Since $\phi_\tau$ vanishes on $\partial\Omega_\tau$, from formula (6.8) follows that

$$0 = \phi_\tau(\boldsymbol{x}_\tau) = \phi_\tau(\boldsymbol{x}_S + \sigma(\boldsymbol{x}_w - \boldsymbol{x}_S)) = \phi_\tau(\boldsymbol{x}_S) + \sigma\nabla\phi_\tau|_K \cdot (\boldsymbol{x}_w - \boldsymbol{x}_S), \tag{6.10}$$

where in the last equality, we use that $\phi_\tau|_K$ is an affine function. Since $\boldsymbol{x}_w$ and each $\boldsymbol{x}_S \in \widehat{S}$ are on opposite sides of $\partial\Omega_\tau$, the directional derivative in the last term of expression (6.10) cannot vanish. We may thus, for each $x_S \in \hat{S}$, solve equation (6.10) for $\sigma$, using that $\phi_\tau(x_S) = \phi(x_S)$, substitute $\sigma$ into expression (6.8), and define a rational function $\hat{X} : \hat{S} \times (0, t)$ by

$$\hat{X}(x_S, \tau) = x_S - \frac{\phi(x_S)}{\nabla \phi_\tau|_K \cdot (x_w - x_S)}(x_w - x_S), \tag{6.11}$$

which satisfies $\hat{X}(\hat{S}, \tau) = \partial\Omega_\tau \cap K$ and $\hat{X}(\hat{S}, (0, t)) = E_t \cap K$. Thus, $\hat{X}$ is surjective. Smoothness follows from that $\hat{X}$ is rational and nonsingular, and injectiveness from the uniqueness of representation (6.8) together with the fact that $\partial\Omega_\tau \cap K$ and $\partial\Omega_{\tau'} \cap K$ do not intersect whenever $\tau \neq \tau'$. □

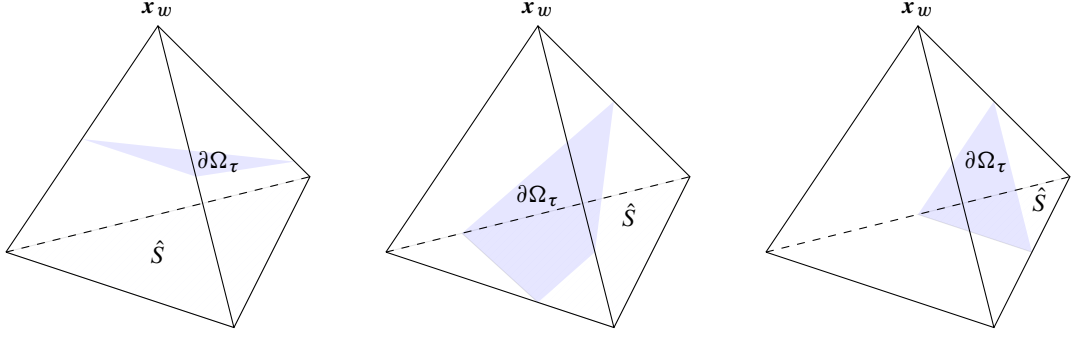

FIGURE 4. Each point $x_\tau$ on the perturbed level-set surface $\partial\Omega_\tau \cap K$, for some element $K \in \mathcal{T}_h$, can be written as $x_\tau = x_S + \sigma(x_w - x_S)$ where $x_w$ is the center node of the basis function $w$ involved in the perturbation, $\sigma \in (0, 1)$, and $x_S \in \hat{S} \subset S$, where $S$ is the mesh surface opposite to $x_w$. The figure illustrates the three different cases in which $\partial\Omega_\tau$ can intersect a tetrahedral mesh element. The domain $\hat{S}$ consists of all points $x_S$ on $S$ such that $\phi(x_S)$ has the opposite sign from $\phi_{t_{\max}}(x_w)$.

Lemma 6.5 provides the necessary tool to parametrize the integral over $E_t$ involved in expression (6.6), starting with the following result.

**Lemma 6.6.** *Consider perturbation (6.2) for some $t \in (0, t_{max}]$, according to assumption 6.3. Let $K \in \mathcal{T}_h$ be an element in the support of basis function $w$. For $f \in C^0(\overline{E_t \cap K})$, it holds that*

$$\int_{E_t \cap K} f \, dV = \int_0^t \int_{\partial\Omega_\tau \cap K} f \frac{w}{|\partial_n \phi_\tau|} \, dS \, d\tau. \tag{6.12}$$

*Proof.* We may again assume that $E_t \cap K$ is nonempty, otherwise the statement is vacuously true. Let $\hat{X}$ be the mapping of lemma 6.5 and $\mathcal{U}^{d-1}$ the domain of a smooth parametrization of the region $\hat{S}$, which is polygonal for $d = 3$ and a line segment for $d = 2$. Then we may devise a smooth parametrization of $E_t \cap K$ with a function $X : \mathcal{U}^{d-1} \times (0, t)$, defined by

$$X(\mathbf{u}, \tau) = \hat{X}(x_S(\mathbf{u}), \tau), \tag{6.13}$$

where $\mathbf{u} = (u, v) \in \mathcal{U}^2$ for $d = 3$ and $\mathbf{u} = u \in \mathcal{U}^1$ for $d = 2$. The parametrization is constructed such that $X(\mathcal{U}^{d-1}, \tau) = \partial\Omega_\tau \cap K$.

The Jacobian determinant of $X$ satisfies

$$\det DX = \frac{\partial X}{\partial \tau} \cdot \mathbf{n}_\tau \sigma_{d-1}; \tag{6.14}$$

where

$$\mathbf{n}_\tau \sigma_2 = \pm \frac{\partial X}{\partial u} \times \frac{\partial X}{\partial v}, \qquad \mathbf{n}_\tau \sigma_1 = \pm Q \frac{\partial X}{\partial u}, \tag{6.15}$$

in which $Q = \begin{pmatrix} 0 & -1 \\ 1 & 0 \end{pmatrix}$; where $\mathbf{n}_\tau$ is the unit normal to $\partial\Omega_\tau \cap K$, outward directed with respect to $\Omega_\tau$ as a function of the parametrization variable $\mathbf{u} \in \mathcal{U}^{d-1}$; where $\sigma_{d-1} > 0$; and where the sign depends on the orientation of the parametrization.

In order to obtain an expression for the right side of expression (6.14), we note that since $X(\mathcal{U}^{d-1}, \tau) = \partial\Omega_\tau \cap K$, by Lemma 6.5, and $\phi_\tau|_{\partial\Omega_\tau \cap K} = 0$, we compose $\phi_\tau$ with $X$ and conclude that on $\mathcal{U}^{d-1} \times (0, t)$,

$$\phi_\tau \circ X = \phi \circ X + \tau w \circ X = 0. \tag{6.16}$$

Differentiating equation (6.16) with respect to $\tau$, we find that

$$\begin{aligned}
\frac{\partial}{\partial \tau}\left(\phi \circ X + \tau w \circ X\right) &= (\nabla\phi \circ X) \cdot \frac{\partial X}{\partial \tau} + w \circ X + \tau(\nabla w \circ X) \cdot \frac{\partial X}{\partial \tau} \\
&= w \circ X + (\nabla\phi_\tau \circ X) \cdot \frac{\partial X}{\partial \tau} = w \circ X + |\nabla\phi_\tau \circ X| \boldsymbol{n}_\tau \cdot \frac{\partial X}{\partial \tau} = 0,
\end{aligned} \tag{6.17}$$

where in the third equality, we use the fact that

$$\boldsymbol{n}_\tau = \frac{\nabla\phi_\tau \circ X}{|\nabla\phi_\tau \circ X|}. \tag{6.18}$$

Combining expressions (6.17) and (6.14), we find that

$$\det DX = -\frac{w \circ X}{|\nabla\phi_\tau \circ X|}\sigma_{d-1}. \tag{6.19}$$

Now we compute using parametrization $X$,

$$\begin{aligned}
\int_{E_t \cap K} f \, dV &= \int_0^t \int_{\mathcal{U}^{d-1}} f \circ X \, |\det DX| \, d\boldsymbol{u} \, d\tau = \int_0^t \int_{\mathcal{U}^{d-1}} f \circ X \frac{w \circ X}{|\nabla\phi_\tau \circ X|}\sigma_{d-1} \, d\boldsymbol{u} \, d\tau \\
&= \int_0^t \int_{\partial\Omega_\tau \cap K} f \frac{w}{|\nabla\phi_\tau|} \, dS \, d\tau = \int_0^t \int_{\partial\Omega_\tau \cap K} f \frac{w}{|\partial_n\phi_\tau|} \, dS \, d\tau,
\end{aligned} \tag{6.20}$$

where expression (6.19) and the fact that that $w \geq 0$ is used in the second equality, and a change of variables is carried out in the third equality. In the last equality, we use that the tangential gradient of $\phi_\tau$ vanishes on $\partial\Omega_\tau \cap K$. $\square$

Now we are ready to work out the directional semiderivative of volume integral $J_1$.

**Theorem 6.7.** *Under perturbation (6.2), the directional semiderivative of volume integral $J_1$ in expression (6.4), for $f \in C^0(\bar{\mathcal{T}}_h)$, satisfies*

$$dJ_1(\phi; w) = \lim_{t \to 0^+} \frac{1}{t}\left(J_1(\phi_t) - J_1(\phi)\right) = -\int_{\partial\Omega} f \frac{w}{|\partial_n\phi|} \, dS. \tag{6.21}$$

(i) *If assumption 6.4 is violated, then $f$ and $\partial_n\phi$ are the limits of these functions from the interior of $\Omega$ whenever these quantities possess jump discontinuities on $\partial\Omega$.*

(ii) *If assumption 6.4 is satisfied, the semiderivatives $t \to 0^-$ and $t \to 0^+$ agree.*

*Proof.* If the support of basis function $w$ does not intersect $\Omega$, formula (6.21) is trivially true. Otherwise, let $t \in (0, t_{\max}]$, according to assumption 6.3, use expression (6.6) and lemma 6.6 to obtain

$$\frac{1}{t}\left(J_1(\phi_t) - J_1(\phi)\right) = -\frac{1}{t}\sum_{K \in \mathcal{T}_h} \int_{E_t \cap K} f \, dV = -\frac{1}{t}\sum_{K \in \mathcal{T}_h} \int_0^t \int_{\partial\Omega_\tau \cap K} f \frac{w}{|\partial_n\phi_\tau|} \, dS \, d\tau. \tag{6.22}$$

Formula (6.21) then follows after passing to the limit in $t$. Conclusion (i) follows from the fact that $\partial\Omega_\tau$ is interior to $\Omega$ for each $\tau \in (0, t)$. As concluded in remark 6.2, an analogous analysis as above can be carried out for $t \leq 0$ which would also lead to formula as (6.21) for $t \to 0^-$, but with $f$ and $\partial_n\phi$ being the limits from the exterior of $\Omega$. Under assumption 6.4, conclusion (ii) then follows since $f$ and $\partial_n\phi$ are then continuous almost everywhere on $\partial\Omega$. $\square$

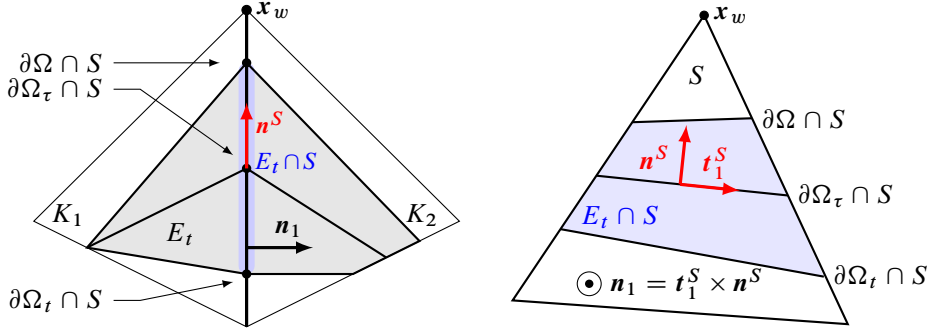

FIGURE 5. The relation between the normal vector $\mathbf{n}_1$ to a mesh surface $\bar{S} = \bar{K}_1 \cap \bar{K}_2$ and the $S$-confined vector $\mathbf{n}^S$ (and $\mathbf{t}_1^S$ for $d = 3$). The central mesh node for the perturbation basis function is marked $\mathbf{x}_w$, here assumed to be outside of $\Omega$.

Left ($d = 2$): $E_t$ is the gray area, $K_1$ is to the left, $K_2$ to the right, and vector $\mathbf{n}_1$ points from $K_1$ to $K_2$. "Surface" $S$ is a line segment; $E_t \cap S$ is the tinted part of $S$; $\partial\Omega_\tau \cap S$ is a point in $S$; and $\mathbf{n}^S$ is the normal to $\partial\Omega_\tau \cap S$, confined in $S$ and orthogonal to $\mathbf{n}_1$.

Right ($d = 3$): $K_1$ and $K_2$ are below and above the plane, respectively, and $\mathbf{n}_1$ points upwards, orthogonal to the plane. Surface $S$ is the triangle; $E_t \cap S$ is the tinted polygonal area; and $\partial\Omega_\tau$ is a line segment with normal $\mathbf{n}^S$ and tangent $\mathbf{t}_1^S$, confined in $S$ and oriented such that $\mathbf{n}_1 = \mathbf{t}_1^S \times \mathbf{n}^S$.

### 6.2 The boundary integral $J_2$

In order to calculate the semiderivative of boundary integral $J_2$ in definition (6.4), the basic idea is to involve the normal field on $\partial\Omega$ and think about the integral as

$$J_2(\phi) = \int_{\partial\Omega(\phi)} f |\mathbf{n}| \, dS. \qquad (6.23)$$

The calculation will use the fact that the normal field on the perturbed boundary $\partial\Omega_t$ can be computed by the formula

$$\mathbf{n}_t = \left.\frac{\nabla \phi_t}{|\nabla \phi_t|}\right|_{\partial\Omega_t}. \qquad (6.24)$$

As previously mentioned, all calculations in this section will be made under assumption 6.4, that is, that the boundary does not intersect any mesh points. A consequence of this assumption is that the right side of expression (6.24) is well defined almost everywhere on $\partial\Omega_t$. However, due to the jump discontinuities of $\nabla \phi_t$ across elements boundaries, the final expression will involve more terms, generated by these jumps, compared to the domain integral case above.

We start by introducing some additional notation that will be used throughout the rest of this section. Figure 5 illustrates the definitions introduced below.

For each face $S \in \mathscr{S}_h$ not residing on $\partial D$, there are two adjacent elements $K_1, K_2 \in \mathscr{T}_h$ such that $\bar{S} = \bar{K}_1 \cap \bar{K}_2$. Denote by $\mathbf{n}_1$ the normal to $S$ directed towards $K_2$ and by $\mathbf{n}_2 = -\mathbf{n}_1$ the normal directed toward $K_1$. Let $t \in (0, t_{\max}]$, and assume that the support of basis function $w$ in perturbation (6.2) has a nonempty intersection with $\Omega$ so that $E_t = \Omega \setminus \overline{\Omega_t}$ is nonempty. Consider adjacent elements $K_1, K_2 \in \mathscr{T}_h$ such that $E_t \cap K_1$ and $E_t \cap K_2$ are both nonempty. For any $\tau \in [0, t]$, we define a unit vector $\mathbf{n}^S$ located in $S$, outward directed from $\Omega_\tau \cap S$; for $d = 3$ we require $\mathbf{n}^S$ to be orthogonal to the edge $\partial\Omega_\tau \cap S$. Note that the outward normal field on $\partial\Omega_\tau$ may have a jump discontinuity at $\partial\Omega_\tau \cap S$ and that vector $\mathbf{n}^S$ is in the span of the normals on $\partial\Omega_\tau \cap K_1$ and $\partial\Omega_\tau \cap K_2$. For $d = 3$, there are also two opposite-directed tangent vectors $\mathbf{t}_1^S, \mathbf{t}_2^S$ in $S$ to $\partial\Omega_\tau \cap S$. Given $\mathbf{n}_k$ and $\mathbf{n}^S$, we choose the orientation of these tangent vectors such that

$$\mathbf{n}_k = \mathbf{t}_k^S \times \mathbf{n}^S \qquad k = 1, 2. \qquad (6.25)$$

Since $\mathbf{n}_2 = -\mathbf{n}_1$, it follows that $\mathbf{t}_2^S = -\mathbf{t}_1^S$. Note that vectors $\mathbf{n}^S, \mathbf{t}_k^S$ are constant on $\partial\Omega_\tau \cap S$ but depend in general on $\tau$, and that $\mathbf{n}_k$ is constant on $S$.

For $d = 2$, $S$ is a line segment, $\mathbf{n}^S$ is parallel to $S$, orthogonal to $\mathbf{n}_k$, and it holds that $\mathbf{n}_2 = Q\mathbf{n}^S$, with $Q = \begin{pmatrix} 0 & -1 \\ 1 & 0 \end{pmatrix}$, and $\mathbf{n}_1 = -\mathbf{n}_2$. However, in order to unify the notation, we will also for $d = 2$ use the notation (6.25), where $\boldsymbol{n}_k$ and $\boldsymbol{n}^S$ are regarded as extended, by zero in the third component, into $\mathbb{R}^3$ using right-handed orientation, and where $\boldsymbol{t}_2^S = \boldsymbol{e}_3$ and $\boldsymbol{t}_1^S = -\boldsymbol{e}_3$.

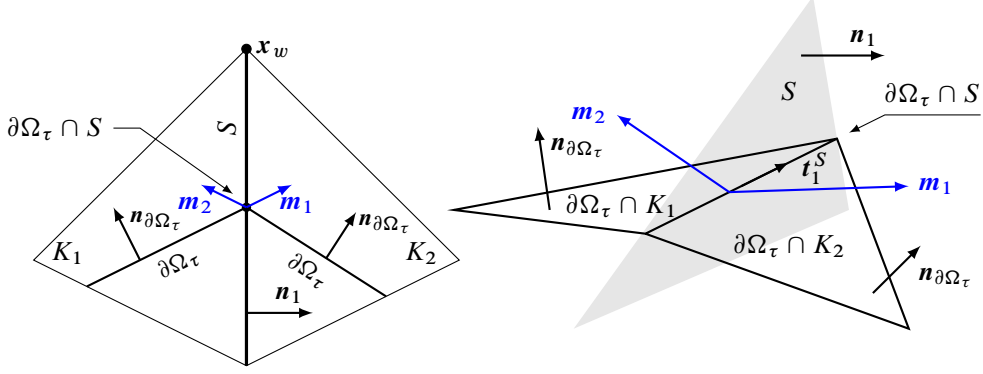

FIGURE 6. Illustration for $d = 2$ (left) and $d = 3$ (right) of the co-normals $\boldsymbol{m}_k$ defined in expression (6.26). These are in the planes of $\partial\Omega_\tau \cap K_k$ and normal to the tangent vector $\boldsymbol{t}_k^S$ along the edge $\partial\Omega_\tau \cap S$.

The last piece of notation concerns co-normals $\boldsymbol{m}_k$ which are normal to the edges of surface patches $\partial\Omega_\tau \cap K_k$, as illustrated in figure 6. Let $\boldsymbol{n}_{\partial\Omega_\tau \cap K_k}$ be the outward-directed, with respect to $\Omega_\tau$, normal field on surface patch $\partial\Omega_\tau \cap K_k$. The co-normals are then defined by

$$\boldsymbol{m}_k = \boldsymbol{t}_k^S \times \boldsymbol{n}_{\partial\Omega_\tau \cap K_k}. \tag{6.26}$$

Note that the co-normals depend on $\tau \in [0, t]$. Given a function $f \in C^0(\bar{\mathcal{T}}_h)$ with jump discontinuities at $S$, vector $\boldsymbol{m}_k$ can be used to specify the limit function $f_k$ of $f$ when approaching $\partial\Omega_\tau \cap S$ from $\partial\Omega_\tau \cap K_k$. That is, for each $\boldsymbol{x} \in \partial\Omega_\tau \cap S$,

$$f_k(\boldsymbol{x}) = \lim_{\epsilon \to 0^+} f(\boldsymbol{x} - \epsilon \boldsymbol{m}_k). \tag{6.27}$$

In the following, we will carry out integration over mesh surfaces that cut through $E_t$, that is, over nonempty $E_t \cap S$ for some $S \in \mathcal{S}_h$. To parametrize $E_t \cap S$, which is an object of dimension $d - 1$, we use an analogous construction as in the proof of Lemma 6.5, which parametrized the $d$-dimensional object $E_t \cap K$. Assumption 6.3, which was used in the proof, has an analogue for each $S \in \mathcal{S}_h$, as remarked in point 1 in the discussion after assumption 6.3. Hence, the following results are proven in the same way as Lemma 6.5.

**Lemma 6.8.** *Consider perturbation (6.2) for some $t \in (0, t_{max}]$, according to assumption 6.3. Let face $S \in \mathcal{S}_h$ be in the support of basis function $w$. Then, for some subface $e \in \mathcal{E}_h$, there exists an $\hat{e} \subset e$ and a diffeomorphism $\widehat{X} : \hat{e} \times (0, t)$ such that $\widehat{X}(\hat{e}, (0, t)) = E_t \cap S$ and $\widehat{X}(\hat{e}, \tau) = \partial\Omega_\tau \cap S$ for each $\tau \in (0, t)$.*

With the help of Lemma 6.8, we may now parametrize integrals over $E_t \cap S$ and arrive at the following result, whose proof is a precise analogue of the proof of Lemma 6.6,

**Lemma 6.9.** *For $S \in \mathcal{S}_h$, $g \in C^0(\overline{E_t \cap S})$, it holds that*

$$\int_{E_t \cap S} g \, dS = \int_0^t \int_{\partial\Omega_\tau \cap S} g \frac{w}{|\partial_{\boldsymbol{n}^S} \phi_\tau|} \, d\gamma \, d\tau. \tag{6.28}$$

Note that $\phi_\tau|_S$ is differentiable in the $\boldsymbol{n}^S$ direction although it is not differentiable in the directions orthogonal to $S$!

The crucial tool for the differentiation of $J_2$ will be the following integration-by-parts formula involving functions $f$ and $\boldsymbol{\Theta}$ containing jump discontinuities along mesh surfaces.

**Lemma 6.10.** *Let $t \in (0, t_{max}]$ according to assumption 6.3 and consider perturbation (6.2). Moreover, let $f \in C^1(\overline{\mathcal{T}_h})$ and $\boldsymbol{\Theta} \in L^\infty(E_t)^d$ such that $\boldsymbol{\Theta}|_{E_t \cap K} \in C^1(\overline{E_t \cap K})^d$ for each $K \in \mathcal{T}_h$. It holds that*

$$\sum_{K \in \mathcal{T}_h} \int_{E_t \cap K} \left[ (\nabla \cdot \boldsymbol{\Theta}) f + \boldsymbol{\Theta} \cdot \nabla f \right] dV$$

$$= \int_{\partial E_t} \boldsymbol{n} \cdot \boldsymbol{\Theta} f \, dS + \sum_{S \in \mathcal{S}_h} \int_0^t \int_{\partial \Omega_\tau \cap S} \boldsymbol{n}^S \cdot [\![ f \boldsymbol{\Theta} \times \boldsymbol{t}^S ]\!] \frac{w}{|\partial_{n^S} \phi_\tau|} \, d\gamma \, d\tau, \quad (6.29)$$

*where*

$$[\![ f \boldsymbol{\Theta} \times \boldsymbol{t}^S ]\!] = f_1 \boldsymbol{\Theta}_1 \times \boldsymbol{t}_1^S + f_2 \boldsymbol{\Theta}_2 \times \boldsymbol{t}_2^S, \quad (6.30)$$

*in which $\boldsymbol{t}_k^S$ is a tangent vector to $\partial \Omega_\tau \cap S$, related, through expression (6.25), to $\boldsymbol{n}^S$ and a given normal vector $\boldsymbol{n}_k$ to the plane $S$. Moreover, $f_k$ and $\boldsymbol{\Theta}_k$ are the limits of $f$ and $\boldsymbol{\Theta}$ on $\partial \Omega_\tau \cap S$ defined by*

$$f_k(\boldsymbol{x}) = \lim_{\epsilon \to 0^+} f(\boldsymbol{x} - \epsilon \boldsymbol{m}_k), \qquad \boldsymbol{\Theta}_k(\boldsymbol{x}) = \lim_{\epsilon \to 0^+} \boldsymbol{\Theta}(\boldsymbol{x} - \epsilon \boldsymbol{m}_k), \qquad \forall \boldsymbol{x} \in \partial \Omega_\tau \cap S, \quad (6.31)$$

*where $\boldsymbol{m}_k$ are the co-normals of expression (6.26).*

*Proof.* We assume the support of basis function $w$ in perturbation (6.2) has a nonempty intersection with $\Omega$ so that $E_t = \Omega \setminus \overline{\Omega}_t$ is nonempty; otherwise the statement is vacuously true. Let $K \in \mathcal{T}_h$, integrate $(\nabla \cdot \boldsymbol{\Theta}) f$ by parts on $E_t \cap K$, and sum over all elements to find

$$\sum_{K \in \mathcal{T}_h} \int_{E_t \cap K} \left[ (\nabla \cdot \boldsymbol{\Theta}) f + \boldsymbol{\Theta} \cdot \nabla f \right] dV$$

$$= \int_{\partial E_t} \boldsymbol{n} \cdot \boldsymbol{\Theta} f \, dS + \sum_{S \in \mathcal{S}_h} \int_{E_t \cap S} \left( f_1 \boldsymbol{n}_1 \cdot \boldsymbol{\Theta}_1 + f_2 \boldsymbol{n}_2 \cdot \boldsymbol{\Theta}_2 \right) dS. \quad (6.32)$$

Note that the last term of formula (6.32) sums solely over mesh surfaces interior to $E_t$ and that this term vanishes if $f$ and the normal component of $\boldsymbol{\Theta}$ are continuous between elements, which will turn expression (6.32) into a "normal" integration-by-parts formula.

By lemma 6.9, the last integral in expression (6.32) satisfies

$$\int_{E_t \cap S} \left( f_1 \boldsymbol{n}_1 \cdot \boldsymbol{\Theta}_1 + f_2 \boldsymbol{n}_2 \cdot \boldsymbol{\Theta}_2 \right) dS = \int_0^t \int_{\partial \Omega_\tau \cap S} \left( f_1 \boldsymbol{n}_1 \cdot \boldsymbol{\Theta}_1 + f_2 \boldsymbol{n}_2 \cdot \boldsymbol{\Theta}_2 \right) \frac{w}{|\partial_{n^S} \phi_\tau|} \, d\gamma \, d\tau. \quad (6.33)$$

By expression (6.25), the properties of the triple product, and definition (6.30), it follows that

$$\begin{aligned} f_1 \boldsymbol{n}_1 \cdot \boldsymbol{\Theta}_1 + f_2 \boldsymbol{n}_2 \cdot \boldsymbol{\Theta}_2 &= f_1 (\boldsymbol{t}_1^S \times \boldsymbol{n}^S) \cdot \boldsymbol{\Theta}_1 + f_2 (\boldsymbol{t}_2^S \times \boldsymbol{n}^S) \cdot \boldsymbol{\Theta}_2 \\ &= \boldsymbol{n}^S \cdot (f_1 \boldsymbol{\Theta}_1 \times \boldsymbol{t}_1^S + f_2 \boldsymbol{\Theta}_2 \times \boldsymbol{t}_2^S) = \boldsymbol{n}^S \cdot [\![ f \boldsymbol{\Theta} \times \boldsymbol{t}^S ]\!] \big|_S. \end{aligned} \quad (6.34)$$

Substituting identity (6.34) into expression (6.33), we find that

$$\int_{E_t \cap S} \left( f_1 \boldsymbol{n}_1 \cdot \boldsymbol{\Theta}_1 + f_2 \boldsymbol{n}_2 \cdot \boldsymbol{\Theta}_2 \right) dS = \int_0^t \int_{\partial \Omega_\tau \cap S} \boldsymbol{n}^S \cdot [\![ f \boldsymbol{\Theta} \times \boldsymbol{t}^S ]\!] \frac{w}{|\partial_{n^S} \phi_\tau|} \, d\gamma \, d\tau, \quad (6.35)$$

which substituted into expression (6.32) yields the result. □

Now we are ready to obtain the final expression for the semiderivative of $J_2$.

**Theorem 6.11.** *Consider perturbation (6.2) according to assumption 6.3 of a domain respecting assumption 6.4. The semiderivative of boundary integral $J_2$ in expression (6.4), for $f \in C^1(\overline{\mathcal{T}_h})$ satisfies*

$$\begin{aligned} dJ_2(\phi, w) &= \lim_{t \to 0^+} \frac{1}{t} \left( J_2(\phi_t) - J_2(\phi) \right) \\ &= -\int_{\partial \Omega} \frac{\partial f}{\partial n} \frac{w}{|\partial_n \phi|} \, dS - \sum_{S \in \mathcal{S}_h} \int_{\partial \Omega \cap S} \boldsymbol{n}^S \cdot [\![ f \boldsymbol{m} ]\!] \frac{w}{|\partial_{n^S} \phi|} \, d\gamma, \end{aligned} \quad (6.36)$$

where $\mathbf{n}^S$ is a unit vector located in $S$, outward-directed from $\Omega$, and orthogonal to $\partial\Omega \cap S$ for $d = 3$. Moreover, $[\![f\mathbf{m}]\!] = f_1\mathbf{m}_1 + f_2\mathbf{m}_2$. Here, $\mathbf{m}_k$ is the conormal for $\tau = 0$ defined in expression (6.26), and $f_k$ is the limit $f$ on $S$ defined by $f_k(\mathbf{x}) = \lim_{\epsilon \to 0+} f(\mathbf{x} - \epsilon\mathbf{m}_k)$ for $\mathbf{x} \in \partial\Omega \cap S$.

*Remark 6.12.* For $d = 2$, the second integral in expression (6.36) should be interpreted as a point evaluation; that is, the expression simply is

$$dJ_2(\phi, w) = -\int_{\partial\Omega} \frac{\partial f}{\partial n} \frac{w}{|\partial_n \phi|} dS - \sum_{S \in \mathscr{S}_h} \mathbf{n}^S \cdot \left( [\![f\mathbf{m}]\!] \frac{w}{|\partial_{n^S} \phi|} \right) \bigg|_{\partial\Omega \cap S}. \tag{6.37}$$

Here $\mathbf{n}^S$ is a unit vector along the face $S$, outward-directed from $\Omega$.

*Proof.* Let $t \in (0, t_{\max}]$ according to assumption 6.3, and assume that the support of basis function $w$ in perturbation (6.2) has a nonempty intersection with $\Omega$ so that $E_t = \Omega \setminus \overline{\Omega}_t$ is nonempty. (Otherwise, formula (6.36) is trivially true.)

Due to assumption 6.4, almost everywhere on $\partial\Omega_\tau$, for $\tau \in [0, t]$, it holds for the normal field on $\partial\Omega_\tau$, outward-directed from $\Omega_\tau$, that

$$\mathbf{n}_\tau = \frac{\nabla \phi_\tau}{|\nabla \phi_\tau|}\bigg|_{\partial\Omega_\tau}. \tag{6.38}$$

Let $K \in \mathscr{T}_h$ have a nonempty intersection with $E_t$. Since the mapping $\widehat{X}$ of lemma 6.5 is an diffeomorphism on $E_t \cap K$ (and thus invertible), we may define the function $\boldsymbol{\psi}_K^{(t)} : E_t \cap K \to \mathbb{R}^N$ through its composition with $\widehat{X}$ such that

$$\boldsymbol{\psi}_K^{(t)}(\widehat{X}(\mathbf{x}_S, \tau)) = \frac{\nabla \phi_\tau(\widehat{X}(\mathbf{x}_S, \tau))}{|\nabla \phi_\tau(\widehat{X}(\mathbf{x}_S, \tau))|} \qquad \forall (\mathbf{x}_S, \tau) \in \widehat{S} \times (0, t), \tag{6.39}$$

We moreover define the assembled function $\boldsymbol{\psi}^{(t)} \in L^\infty(E_t)^N$ such that $\boldsymbol{\psi}^{(t)}|_{E_t \cap K} = \boldsymbol{\psi}_K^{(t)}$ for each $K \in \mathscr{T}_h$ having nonempty intersection with $E_t$. Function $\boldsymbol{\psi}^{(t)}$ is constructed, using property (6.38), so that, for $K$ having nonempty intersection with $E_t$,

$$\boldsymbol{\psi}^{(t)}|_{\partial\Omega_\tau \cap K} = \mathbf{n}_\tau|_{\partial\Omega_\tau \cap K} \qquad \forall \tau \in (0, t). \tag{6.40}$$

By the smoothness of mapping $\widehat{X}$, function $\boldsymbol{\psi}^{(t)}$ fulfills the conditions for $\boldsymbol{\Theta}$ in lemma 6.10, which implies that for $f \in C^1(\overline{\mathscr{T}}_h)$,

$$\sum_{K \in \mathscr{T}_h} \int_{E_t \cap K} [(\nabla \cdot \boldsymbol{\psi}^{(t)}) f + \boldsymbol{\psi}^{(t)} \cdot \nabla f] \, dV$$
$$= \int_{\partial E_t} \mathbf{n} \cdot \boldsymbol{\psi}^{(t)} f \, dS + \sum_{S \in \mathscr{S}_h} \int_0^t \int_{\partial\Omega_\tau \cap S} \mathbf{n}^S \cdot [\![f\boldsymbol{\psi}^{(t)} \times \mathbf{t}^S]\!] \frac{w}{|\partial_{n^S} \phi_\tau|} \, d\gamma \, d\tau. \tag{6.41}$$

Using lemma 6.6 and property (6.40), nonvanishing element contributions to the left side of formula (6.41) satisfy

$$\int_{E_t \cap K} [(\nabla \cdot \boldsymbol{\psi}^{(t)}) f + \boldsymbol{\psi}^{(t)} \cdot \nabla f] \, dV = \int_0^t \int_{\partial\Omega_\tau \cap K} (\underbrace{\nabla \cdot \mathbf{n}_\tau}_{=\kappa_\tau} f + \mathbf{n}_\tau \cdot \nabla f) \frac{w}{|\partial_n \phi_\tau|} d\gamma \, d\tau$$
$$= \int_0^t \int_{\partial\Omega_\tau \cap K} \mathbf{n}_\tau \cdot \nabla f \frac{w}{|\partial_n \phi_\tau|} d\gamma \, d\tau, \tag{6.42}$$

where $\kappa_\tau$ is the summed curvature of $\partial\Omega_\tau \cap K$, which vanishes since $\partial\Omega_\tau \cap K$ is planar.

Regarding the first term on the right side of formula (6.41), note that $\partial E_t$ consists of two parts, corresponding to the boundaries of $\partial\Omega$ and $\partial\Omega_t$, respectively (cf. figure 3). Moreover, the integrand of this term satisfies

$$\boldsymbol{n}\cdot\boldsymbol{\psi}^{(t)}f\big|_{\partial E_t} = \begin{cases} f & \text{on } \partial\Omega\cap E_t, \\ -f & \text{on } \partial\Omega_t\cap E_t, \end{cases} \qquad (6.43)$$

where the minus sign is due to that $\boldsymbol{n}$ here denotes the outward-directed normal from $E_t$, whereas $\boldsymbol{\psi}^{(t)}$, due to property (6.40), is directed outward from $\Omega_\tau$; those directions are the same on $\partial\Omega$ but opposite on $\partial\Omega_t$.

Expression (6.43) implies that the first term on the right side of formula (6.41) can be written

$$\int_{\partial E_t} \boldsymbol{n}\cdot\boldsymbol{\psi}_t f\, dS = \int_{\partial\Omega\cap E_t} f\, dS - \int_{\partial\Omega_t\cap E_t} f\, dS = \int_{\partial\Omega} f\, dS - \int_{\partial\Omega_t} f\, dS, \qquad (6.44)$$

where the last equality follows from that $\partial\Omega = \partial\Omega_t$ outside of $E_t$.

Recalling jump notation (6.30) and property (6.40), the triple product in the last integral of formula (6.41) can be written

$$\begin{aligned}
\boldsymbol{n}^S \cdot [\![f\boldsymbol{\psi}^{(t)}\times\boldsymbol{t}^S]\!]\big|_{\partial\Omega_\tau\cap S} &= \boldsymbol{n}^S\cdot\left(f_1\boldsymbol{\psi}_1^{(t)}\times\boldsymbol{t}_1^S + f_2\boldsymbol{\psi}_2^{(t)}\times\boldsymbol{t}_2^S\right)\Big|_{\partial\Omega_\tau\cap S} \\
&= \boldsymbol{n}^S\cdot\left(f_1\boldsymbol{n}_{\tau,1}\times\boldsymbol{t}_1^S + f_2\boldsymbol{n}_{\tau,2}\times\boldsymbol{t}_2^S\right)\Big|_{\partial\Omega_\tau\cap S} \\
&= \boldsymbol{n}^S\cdot[\![f\boldsymbol{n}_\tau\times\boldsymbol{t}^S]\!]\big|_{\partial\Omega_\tau\cap S} = -\boldsymbol{n}^S\cdot[\![f\boldsymbol{m}]\!]\big|_{\partial\Omega_\tau\cap S},
\end{aligned} \qquad (6.45)$$

where $\boldsymbol{n}_{\tau,k} = \boldsymbol{n}_\tau\big|_{\partial\Omega_\tau\cap K_k}$ and $[\![f\boldsymbol{m}]\!] = f_1\boldsymbol{m}_1 + f_2\boldsymbol{m}_2$, in which $\boldsymbol{m}_k$ denotes the co-normals (6.26), illustrated in figure 6, and where $f_k$ is defined in expression (6.31).

Substituting expressions (6.42), (6.44), and (6.45) into formula (6.41) and rearranging, we find that

$$\begin{aligned}
\int_{\partial\Omega} f\, dS &- \int_{\partial\Omega_t} f\, dS \\
&= \int_0^t \int_{\partial\Omega_\tau} \boldsymbol{n}_\tau\cdot\nabla f \frac{w}{|\partial_n\phi_\tau|}\, d\tau\, dS + \sum_{S\in\mathscr{S}_h} \int_0^t \int_{\partial\Omega_\tau\cap S} \boldsymbol{n}^S\cdot[\![f\boldsymbol{m}]\!]\frac{w}{|\partial_{n^S}\phi_\tau|}\, d\tau\, dS.
\end{aligned} \qquad (6.46)$$

Finally, dividing by $t$ and passing to the limit, we find

$$\begin{aligned}
dJ_2(\phi,w) &= \lim_{t\to 0^+}\frac{1}{t}\bigl(J_2(\phi_t) - J_2(\phi)\bigr) = \lim_{t\to 0^+}\frac{1}{t}\left(\int_{\partial\Omega_t} f\, dS - \int_{\partial\Omega} f\, dS\right) \\
&= -\int_{\partial\Omega}\frac{\partial f}{\partial n}\frac{w}{|\partial_n\phi|}\, dS - \sum_{S\in\mathscr{S}_h}\int_{\partial\Omega\cap S} \boldsymbol{n}^S\cdot[\![f\boldsymbol{m}]\!]\frac{w}{|\partial_{n^S}\phi|}\, dS,
\end{aligned} \qquad (6.47)$$

which is the claim. $\square$

## Acknowledgments

Thanks to Linus Hägg for a careful reading of the manuscript and for providing detailed and useful feedback. Thanks also to Stephan Schmidt and Eddie Wadbro for providing an alternative reasoning to prove the expressions of theorem 6.11. I am also grateful to the anonymous reviewer for a close reading of the manuscript, which detected many small issues to be fixed, and for suggesting the use of a better model problem to demonstrate also the use of the boundary-integral formulas. Support for this work was provided by the Swedish Research Council, grant 2018-03546.